\documentclass[12pt,3p]{elsarticle}
\usepackage{amsfonts,amsmath,amssymb,amsthm,framed}
\usepackage[utf8]{inputenc}
\usepackage{spverbatim}
\usepackage{multicol}
\usepackage{blkarray}
\usepackage{mathtools} 
\usepackage{float}
\usepackage{algpseudocode}
\usepackage[toc]{appendix}

\usepackage{hyperref}
\numberwithin{equation}{section}
\usepackage{xcolor}
\usepackage{tikz}
\usepackage{colortbl}
\usepackage[normalem]{ulem}
\usepackage{sectsty}
\definecolor{astral}{RGB}{46,116,181}
\subsectionfont{\color{astral}}
\sectionfont{\color{astral}}
\usepackage{colortbl}
\newtheorem{theorem}{Theorem}[section]
\newtheorem{lemma}[theorem]{Lemma}

\usepackage{hyperref}

\usepackage{cancel}
\usepackage{graphicx}
\usepackage{enumitem}
\usepackage{fullpage}
\definecolor{warmblack}{rgb}{0.0, 0.26, 0.26}
\usepackage{bm}
\usepackage{tikz}
\usetikzlibrary{ decorations.markings}
\usetikzlibrary{arrows.meta}
\usepackage{array}
\usepackage{booktabs}
\usepackage{tabularx}

\usepackage{amsbsy}
\usepackage{accents}
\newlength{\dhatheight}
    
\newenvironment{frcseries}{\fontfamily{frc}\selectfont}{}

\begin{document}

\begin{frontmatter}

\title{ \textcolor{warmblack}{\bf Network-Driven Global Stability Analysis: SVIRS Epidemic Model 
}}

\cortext[cor1]{Corresponding author}

\author[label1]{Madhab Barman\corref{cor1}}
\address[label1]{Department of Mathematics, Indian Institute of Information Technology Design and Manufacturing Kancheepuram, Chennai-600127, India}
\ead{madhabbrmn@gmail.com}
\author[label1]{Nachiketa Mishra}
\begin{abstract}
An epidemic Susceptible-Vaccinated-Infected-Removed-Susceptible (SVIRS) model is presented on a weighted-undirected network with graph Laplacian diffusion. Disease-free equilibrium always exists while the existence and uniqueness of endemic equilibrium have been shown.  When the basic reproduction number is below unity, the disease-free equilibrium is asymptotically globally stable. The endemic equilibrium is asymptotically globally stable if the basic reproduction number is above unity. Numerical analysis is illustrated with a road graph of the state of Minnesota. The effect of all important model parameters has been discussed. 
\end{abstract}
\end{frontmatter}
\section{Introduction} 
\noindent Any transmission disease with a slow rate of recovery or slow prevention can cause an outbreak in the population. It has a massive influence on human beings. Because of this, every decade millions of citizens lose their life. One of the most cost-effective strategies to mitigate and control the spread of infectious diseases is the vaccine. Vaccines are biological preparations that help the immune system for developing protection from disease. Mathematical compartmental models such as SIR, SIS, SEIR, etc., are often used to analyze how to control the spreading of infectious disease, see \cite{SIR:1927, Myli:1995, book:brauer, AMS:2008}. Many authors have introduced various epidemic models by incorporating a vaccine compartment to investigate the effect of vaccination, see \cite{MB_network,Vac:2022, Mahmood:2018a, Mahmood:2018b, Feng:2011, Liu:2008, sunita:2008, AIP:2021, Li:2011}.
Another crucial component to comprehend the spatial spread of an infectious disease is population mobility. Over the years, several models have been proposed to study the spread of infectious diseases on a complex network. More details about the epidemic process in complex networks can be found at \cite{rev1, rev2}. The Laplacian diffusion has been used recently to understand spatial disease transmission dynamics by researchers \cite{Bo2019,  sir_networ21, network_europe20, network_us20, Ding21, Brock18}. It is also used to characterize population mobility through a complex network. So, adopting Laplacian diffusion is vital for disease transmission dynamics. In recent years, the authors Tian, Barman, Zhang, and Ruan \cite{global:2020,MB_bifur,tian2020,SEIR:2022} have studied different types of networked SIR and SEIR epidemic models. 
\\[2ex]
\noindent Consider $G:=\langle \mathcal{V} , \mathcal{E}, \mathcal{W} \rangle$ is a connected weighted-undirected graph, where $\mathcal{V}$ is the set of nodes or vertices with $|\mathcal{V}| = n$, $\mathcal{E}$ is the set of edges  and $\mathcal{W}$ is the set of weights corresponding $\mathcal{E}$. Taking the graph as a spatial domain, we define the Laplacian operator, $\Delta$, acting on a function $F$  from continuous space to a finite weighted-undirected graph without self loop, see \cite{global:2020,MB_bifur,tian2020}, as below:
\begin{align}\label{delta}
    \Delta F(x) = \sum_{{y \in \mathcal{V}, y \sim x}}w(x, y) \left[ F(y) - F(x)\right],
\end{align}
where $y \sim x$ describes node $y$ is adjacent to node $x$, and $F$ is a function such that $F : \mathcal{V} \to \mathbb{R}$. Consider $w(x, y) > 0$ be the weight corresponding to an edge whose adjacent nodes are $x$ and $y$ such that $w(x, y) = w(y, x) $. Let $S(x, t), V(x,t), I(x, t)$ and $R(x,t)$ are the total number of Susceptible, Vaccinated, Infected, and Removal individuals at node $x \in \mathcal{V}$ and time $t \in [0, +\infty)$. Let $S,V,I,R : \mathcal{V}\times [0, +\infty) \to \mathbb{R}$, then SVIRS model on a Laplacian network can be defined as: 
\begin{equation}
\left\lbrace
    \begin{aligned}\label{model}
    \dfrac{dS}{dt}(x, t) - \epsilon \Delta S(x, t) &= \mu(1-r)N -(\mu+p)S(x, t) - \beta S(x, t) I(x, t)\\
    &~~~+ \xi V(x, t) + \eta R(x, t),\\
    \dfrac{dV}{dt}(x, t) - \epsilon \Delta V(x, t) &=\mu r N + pS(x, t) - \beta \sigma V(x, t) I(x, t) - (\mu+\xi) V(x, t),\\
    \dfrac{dI}{dt}(x, t) - \epsilon \Delta I(x, t) &= \beta S(x, t) I(x, t) +  \beta \sigma V(x, t) I(x, t) - (\mu+\gamma) I(x, t),\\
    \dfrac{dR}{dt}(x, t) - \epsilon \Delta R(x, t) &= \gamma I(x, t) - (\mu+\eta) R(x, t),
    \end{aligned}
    \right.
\end{equation}
 for all $(x, t) \in \mathcal{V}\times [0, +\infty)$,  subject to the initial conditions: 
\begin{equation}\label{init}
 \begin{aligned}
    \Big\lbrace &S(x, 0) = S_0(x), V(x, 0) = V_0(x), I(x, 0) = I_0(x), R(x, 0) = R_0(x)\\
    &\mbox{ and } S_0(x) + V_0(x) + I_0(x) + R_0(x) = N\,\mbox{ for } x \in \mathcal{V}  \Big\rbrace \geq 0.
\end{aligned}   
\end{equation}
All the parameters used in this model are assumed to be non-negative. The descriptions of each parameter are depicted in Table-\ref{tab:parameters}. The model (\ref{model}) without network is borrowed from the article \cite{Liu:2008} by the author Xianning Liu et al.. In the model (\ref{model}), we incorporated three important cases (i)  post-recovery waning rate, $\eta$, refers to the rate at which immunity acquired after recovering from an infectious disease or through vaccination gradually diminishes over time. As time passes, individuals may become more susceptible to reinfection or less protected against the disease, (ii) vaccination for the newborns (in our model, there is an inflow of newborns at a rate of $\mu N$. The parameter $r$ represents the fraction of vaccinated newborns, denoted as $r \mu N$, while $1-r$ denotes the remaining fraction, which comprises the unvaccinated newborns, i.e., $(1-r)\mu N$), and (iii) all individual movement through a weighted-undirected network. The parameter $\epsilon$ is often called migration parameter or diffusion parameter or population mobility parameter,    see \cite{sir_networ21, network_europe20, network_us20}. This parameter helps us to control the population mobility in the network by choosing a suitable value for $\epsilon \in(0, 1]$. When $\epsilon$ is close to $0$, we say there is a strict restriction on the movement of individuals from one location to another location. When $\epsilon = 1$, there is no restriction on the movement of individuals, i.e., individuals can move freely from one place to another. When $\epsilon$ takes a moderate value between 0 and 1, a portion of the population is permitted to move from one location to another.\\[2ex]
\noindent From the assumptions of the model \eqref{model}, one may obtain  $S(x, t) + V(x, t) + I(x, t) + R(x, t)  = N$ (constant) for all $(x, t) \in \mathcal{V} \times [0, +\infty)$. Applying the Lemma 2.6 of \cite{SEIR:2022}, solutions of the model \eqref{model} with the initial condition \eqref{init} are nonnegative on $\mathcal{V} \times [0, T],\, T> 0$. Thus, it is obvious that solutions are bounded and non-negative.\\[2ex]
\noindent The paper is structured as follows: in Section \ref{equilibria}, we find the equilibria of the proposed model \eqref{model} and their existence and uniqueness. Global stability analysis of the equilibria has been presented in Section \ref{global} with the help of the Lyapunov function. Section \ref{numerical} is devoted to numerical results and their discussions. Finally, Section \ref{conclusion} gives a brief conclusion of this paper.

renewcommand{\arraystretch}{1.5}
\begin{table}[h!]
    \centering
    \caption{Descriptions of all the model parameters used in the model (\ref{model})}
{ \begin{tabular}{m{1.5cm} c m{6cm} c m{2.5cm} }
    \toprule
\hline
   {\textbf{Parameters}} & \vline  & { \textbf{Descriptions}} & \vline  & \textbf{Units/Ranges}\\
    \hline
    \hline
         $N$ &:&Total population for each node  &:& $N > 0$\\
        $r$ &:& The fraction of new individuals who are vaccinated &:& $0 < r < 1$\\
        $\mu$ &:& The natural birth and death rate &:& day$^{-1}$\\
        $p $ &:& The proportion of vaccination &:& $0 \leq p \leq 1 $\\
        $\beta$ &:& The rate of disease transmission &:& day$^{-1}$\\
        $\xi$ &:& The rate of losing immunity who are vaccinated  &:& day$^{-1}$\\
         $\eta$ &:&  Post-recovery immunity waning rate  &:& day$^{-1}$\\
         $\sigma $ &:& The effectiveness of vaccination  &:& $0 \leq \sigma \leq 1$\\
         $\gamma $ &:& The removal rate  &:& day$^{-1}$\\
          $\epsilon $ &:& The migration rate or population mobility rate &:& $0<\epsilon \leq 1$\\
         \hline
         \bottomrule
    \end{tabular}}
    \label{tab:parameters}
\end{table}
\renewcommand{\arraystretch}{1}
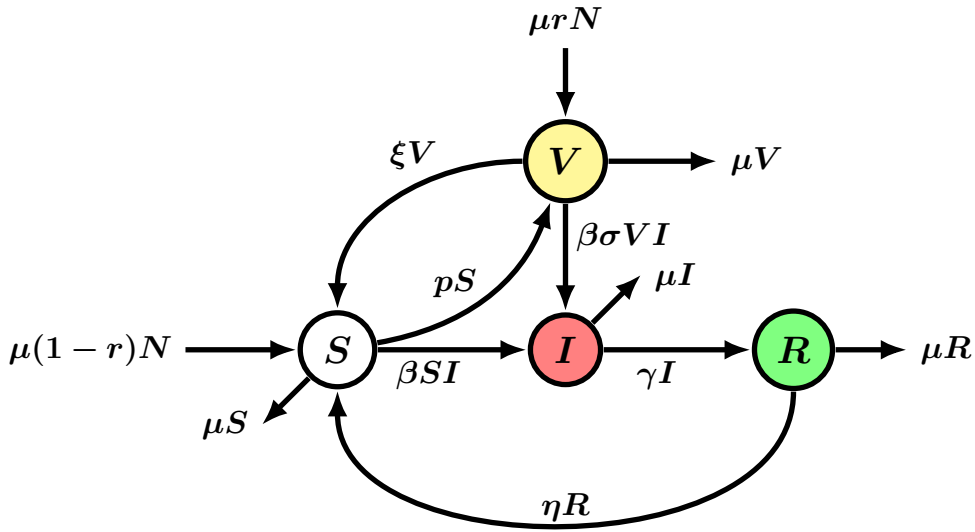
\begin{figure}[ht]
    \centering
      \begin{tikzpicture}
  [decoration={markings, 
    mark= at position 0.25 with {\arrow{stealth}},
    mark= at position 2cm with {\arrow{stealth}}}
] 
\path (0,0) node[circle, draw, fill=white!50, scale=1.2, line width=2pt] (S){ \bm{$S$}}
(3,0) node[circle, draw, fill=red!50, scale=1.2, line width=2pt]  (I){\bm{$I$}}
(6,0) node[circle, draw, fill=green!50, scale=1.2, line width=2pt]  (R){\bm{$R$}}
(3,2.5) node[circle, draw, fill=yellow!50, scale=1.2, line width=2pt]  (V){\bm{$V$}};
\draw[arrows = {-Latex[line width=2pt, fill=black, length=10pt]}, line width=2pt] (S)--(I);
\draw[arrows = {-Latex[line width=2pt, fill=black, length=10pt]}, line width=2pt] (I)--(R);
\draw[arrows = {-Latex[line width=2pt, fill=black, length=10pt]}, line width=2pt] (V)--(I);
\draw[arrows = {-Latex[line width=2pt, fill=black, length=10pt]}, line width=2pt] (-2,0) node[left] {\bm{$\mu(1-r)N$}} --(S);
\draw[arrows = {-Latex[line width=2pt, fill=black, length=10pt]}, line width=2pt] (S)--(-1,-1)node[left]{\bm{$\mu S$}};
\draw[arrows = {-Latex[line width=2pt, fill=black, length=10pt]}, line width=2pt] (I)--(4,1)node[right, yshift=-1pt]{\bm{$\mu I$}};
\draw[arrows = {-Latex[line width=2pt, fill=black, length=10pt]}, line width=2pt] (R)--(7.5,0)node[right]{\bm{$\mu R$}};
\draw[arrows = {-Latex[line width=2pt, fill=black, length=10pt]}, line width=2pt] (V)--(5,2.5) node[right]{\bm{$\mu V$}};
\draw[arrows = {-Latex[line width=2pt, fill=black, length=10pt]}, line width=2pt] (3,4) node[above]{\bm{$\mu r N$}}--(V);
\draw[arrows = {-Latex[line width=2pt, fill=black, length=10pt]}, line width=2pt] (V) to[out=180,in=90](S);
\draw[arrows = {-Latex[line width=2pt, fill=black, length=10pt]}, line width=2pt] (S) to[out=10,in=-110](V);
\draw[arrows = {-Latex[line width=2pt, fill=black, length=10pt]}, line width=2pt] (R) to[out=-90,in=-90](S);
\node[below] at (1.2, 0){\bm{$\beta S I$}};
\node[right] at (3, 1.5){\bm{$\beta \sigma V I$}};
\node[below] at (4.2, 0){\bm{$\gamma I$}};
\node[above] at (1, 2.25){\bm{$\xi V$}};
\node[above] at (1.55, 0.55){\bm{$p S$}};
\node[below] at (3, -1.75){\bm{$\eta R$}};
\end{tikzpicture}
    \caption{Flow chart of the model (\ref{model}) at node $x \in \mathcal{V}$.}
    \label{fig:flowchart}
\end{figure}
section{Equilibria}\label{equilibria}
\noindent For equilibria of the above system (\ref{model}), we set $\dfrac{dS}{dt} = 0, \dfrac{dV}{dt} = 0, \dfrac{dI}{dt} = 0$ and $\dfrac{dR}{dt} = 0$. Considering all the equilibrium points are uniformly constant throughout the network, then $\Delta u = 0$ at the equilibrium point, where $u = S, V, I, R$.  The third equation of system (\ref{model}) gives
\begin{align}
     I\big(\beta (S  +  \sigma V)  - (\mu+\gamma)\big) = 0.
\end{align}
Therefore, either $I = 0$ or $\left(\beta (S  +  \sigma V)  - (\mu+\gamma)\right) = 0$. When $I = 0$, we have the following system
\begin{equation}\label{disease_free_sys}
\left\lbrace
    \begin{aligned}
      0 &= \mu(1-r)N -(\mu+p)S  + \xi V + \eta R,\\
    0 &=\mu r N + pS  - (\mu+\xi) V,\\
    0 &=  - (\mu+\eta) R.
    \end{aligned}
    \right.
\end{equation}
We denote disease-free equilibrium point is $\textbf{D}$ such that $\textbf{D} (S_0, V_0, I_0, R_0)$. Then solving the system (\ref{disease_free_sys}), we acquire  
\begin{align*}
    S_0 = \dfrac{ \xi + \mu(1-r)}{\mu+p+\xi}N,\quad  V_0 = \dfrac{p+\mu r}{\mu+p+\xi}N,\quad  I_0=0,\quad  R_0 = 0.
\end{align*}
\noindent \emph{Basic reproduction number}: { Rewrite the 3rd equation of \eqref{model} into the following ordinary differential equations as follows:
\begin{align*}
    \dfrac{d I_i}{dt}(t) = \beta S_i(t)I_i(t) + \beta\sigma V_i(t) I_i(t) - (\mu+\gamma) I_i(t) + \epsilon \sum_{j=1}^n L_{ij} I_j(t),
\end{align*}
where $I_i(t)$ is the infected individual at node $i$ and time $t$;  $i = 1,2,\dots,n$. Then, at  disease-free point $\textbf{D} $:
\begin{align*}
    \dfrac{d I_i}{dt}(t) = (\beta S_0 + \beta\sigma V_0) I_i(t) - (\mu+\gamma) I_i(t) + \epsilon \sum_{j=1}^n L_{ij} I_j(t),
\end{align*}
The above can be rewritten into vector form as follows:}
\begin{align*}
    \dfrac{d}{dt}\textbf{I} = \Big( (\beta S_0 + \sigma \beta V_0)eye(n) - \big((\mu + \gamma)eye(n)- \epsilon L\big)   \Big) \textbf{I},
\end{align*}
where $eye(n)$ is the identity matrix of size $n$, and $L$ is graph Laplacian matrix associated with the graph $G$ defined by
\begin{align*}
    L_{ij} = \begin{cases}
        ~~~A_{ij}w_{ij}, \qquad~~~~~  j\neq i,\\
        -\sum_{j =1}^n A_{ij} w_{ij}, \quad j=i,
    \end{cases}
\end{align*}
with $A$ is the adjacent matrix of the graph $G$, {and $\textbf{I} = (I_1, I_2, \dots, I_n)^\texttt{T}$; $\texttt{T}$ describes the transpose.} According to the {\it next-generation matrix method} \cite{next_gen},
 we have the  transmissions matrix $\mathbb{T}$, and transitions matrix $\Sigma$:
 \begin{align*}
     \mathbb{T} = (\beta S_0 + \sigma \beta V_0)eye(n), \mbox{ and } \Sigma = (\mu + \gamma)eye(n) -\epsilon L.
 \end{align*}
 Therefore, the basic reproduction number is obtained by
 \begin{align}\label{r0_m}
     \mathcal{R}_0 = \rho\Big(\mathbb{T} \Sigma^{-1} ) = \rho((\beta S_0 + \sigma \beta V_0)\big((\mu + \gamma)eye(n) - \epsilon L \big)^{-1}\Big).
 \end{align}
 Let $\lambda_i, i = 1,2,\dots, n$ be the eigenvalues of $-L$, then $0=\lambda_1 < \lambda_2 \leq \lambda_3\leq \cdots \leq \lambda_n$. Our goal is to find the maximum eigenvalue of the matrix $(\beta S_0 + \sigma \beta V_0)\big((\mu + \gamma)eye(n) - \epsilon L \big)^{-1}$. The smallest eigenvalue of the $(\mu + \gamma)eye(n) -\epsilon L$ is $(\mu +\gamma)$. 
Thus, the basic reproduction number of the entire model is given by 
\begin{align}\label{reprod}
    \mathcal{R}_0 = \frac{\beta (S_0 + \sigma V_0)}{\mu+\gamma} = \frac{\beta \big( ( \xi+(1-r)\mu) + \sigma(p+\mu r)  \big)N}{(\gamma+\mu)(\mu+p + \xi)}.
\end{align}
The endemic equilibrium point, denote it as $\textbf{E} (S_*, V_*, I_*, R_*)$, can be obtained by solving the following system 
\begin{equation}
\left\lbrace
    \begin{aligned}\label{endemic_sys}
    0 &= \mu(1-r)N -(\mu+p)S - \beta S I + \xi V + \eta R,\\
    0 &=\mu r N + pS - \beta \sigma V I - (\mu+\xi) V,\\
    0 &= \beta S +  \beta \sigma V - (\mu+\gamma), \\
    0 &= \gamma I - (\mu+\eta) R.
    \end{aligned}
    \right.
\end{equation}
From the last equation of the system (\ref{endemic_sys}), $R = \gamma I/(\mu+\eta)$. Next, we solve $S$ and $V$ in terms of $I$ using first two equations of the system (\ref{endemic_sys}), 
\begin{equation}
\left\lbrace
    \begin{aligned}
      0 &=  -(\mu+p +\beta I)S  + \xi V + \eta R + \mu(1-r)N, \\
      0 &=   pS - (\beta \sigma I  + \mu+\xi) V + \mu r N. 
    \end{aligned}
    \right.
\end{equation}
A few direct calculations yield
\begin{equation}\label{SV}
\left\lbrace
\begin{aligned}
    S &= \dfrac{\mu N\big(\xi +(1-r)(\mu+\beta \sigma I)  \big) + (\mu + \xi+\beta \sigma I) \eta}{(\mu+p+\beta)(\mu+\xi+\beta \sigma I) -\xi p},\\
     V &= \dfrac{\mu N\big( p+r(\mu+\beta I) \big) + p \eta R}{(\mu+p+\beta)(\mu+\xi+\beta \sigma I) -\xi p}.
\end{aligned}
\right.
\end{equation}
Substituting (\ref{SV}) into third equation of the system (\ref{endemic_sys}) then substitute $R$, we acquire
\begin{align}
    aI^2+bI+c = 0,
\end{align}
where
\begin{equation}\label{abc}
\left\lbrace
    \begin{aligned}
        a &= \frac{\beta \sigma \mu }{\mu+\eta} \left( \mu+\gamma+\eta\right),\\
        b&= (\mu+\gamma)\left[ \sigma(\mu+p)+(\mu+\sigma)\right] - \mu \sigma\beta N - (\mu+\xi+\sigma p)\frac{\gamma \eta}{\mu+\eta},\\
        c &= \frac{\mu+\gamma}{\beta}\big( (\mu+p)(\mu+\xi)-\xi p \big) - \mu N \big(\xi +(1-r)\mu  \big) - \sigma \mu N \big(p +r\mu  \big),\\
        &= \frac{\mu(\mu+\gamma)(\mu+\xi+p)}{\beta} (1 - \mathcal{R}_0).
    \end{aligned}
    \right.
\end{equation}
The roots of the quadratic equation in $I$ are $(-b\pm \sqrt{b^2-4ac})/2a$. From (\ref{abc}), it is clear that $a$ is positive ($a>0$) and $c$ is negative ($c<0$) when $\mathcal{R}_0 > 1$. Therefore, $\sqrt{b^2-4ac}$ will be a positive real number, since $4ac < 0$. The quadratic equation produces one negative root and another positive root because $\sqrt{b^2-4ac} > b$. Again since the negative root does not lie in the feasible region of $I$, the equation will have only one unique positive root when $\mathcal{R}_0 > 1$. Hence the endemic equilibrium exists if $\mathcal{R}_0 > 1$.\\[2ex]
\noindent The global dynamics of the model (\ref{model}) are decided by the reproduction number $\mathcal{R}_0$. In the following section, we present the global stability analysis for both equilibria: disease-free equilibrium (\textbf{D}) and endemic equilibrium (\textbf{E}). 
\section{Global Stability}\label{global}
Before proving theorems of the global analysis, we reduce our model (\ref{model}) by
substituting $R = N - S - V - I$,
\begin{equation}
\left\lbrace
    \begin{aligned}\label{model_reduce}
    \dfrac{dS}{dt} - \epsilon \Delta S &= \mu(1-r)N -(\mu+p)S - \beta S I + \xi V + \eta (N - S - V - I),\\
    \dfrac{dV}{dt} - \epsilon \Delta V &=\mu r N + pS - \beta \sigma V I - (\mu+\xi) V,\\
    \dfrac{dI}{dt} - \epsilon \Delta I &= \beta S I +  \beta \sigma V I - (\mu+\gamma) I.
    \end{aligned}
    \right.
\end{equation}
The above model (\ref{model_reduce}) and the proposed model (\ref{model}) have the same behavior. Since no equation contains the term $R$, the model (\ref{model_reduce}) gives global analysis for $(S, V, I)$.

\begin{lemma}\label{Green}
For any functions $F, G: \mathcal{V}\times [0, \infty) \to \mathbb{R}$ and for fix $t$, we have
\begin{align*}
    \sum_{x \in \mathcal{V}} F(x, t) \Delta G(x, t) = - \frac{1}{2} \sum_{x, y \in \mathcal{V}} w(x, y) \left[ F(y, t) - F(x, t)\right]\left[ G(y, t) - G(x, t)\right].
\end{align*}
\end{lemma}
\noindent\textbf{Proof.} Begin with by dividing the left-hand-side, $\sum\limits_{x \in \mathcal{V}} F(x, t) \Delta G(x, t)$, into two halves as follows:
\begin{align*}
     \sum_{x \in \mathcal{V}} F(x, t) \Delta G(x, t) &= \frac{1}{2} \sum_{x \in \mathcal{V}} F(x, t) \Delta G(x, t) + \frac{1}{2}\sum_{y \in \mathcal{V}} F(y, t) \Delta G(y, t).
\end{align*}
Applying the formula (\ref{delta}), we obtain
\begin{align*}
     \sum_{x \in \mathcal{V}} F(x, t) \Delta F(x, t) &= \frac{1}{2} \sum_{x \in \mathcal{V}} F(x, t) \sum_{{y \in \mathcal{V}, y \sim x}}w(x, y) \left[ G(y, t) - G(x, t)\right] \\
     &~~~+ \frac{1}{2}\sum_{y \in \mathcal{V}} F(y, t) \sum_{{x \in \mathcal{V}, x \sim y}}w(y, x) \left[ G(x, t) - G(y, t)\right].
\end{align*}
Using the symmetry $w(x, y) = w(y ,x)$, 

\begin{align*}
     \sum_{x \in \mathcal{V}} F(x, t) \Delta G(x, t) &= \frac{1}{2} \sum_{x \in \mathcal{V}}  \sum_{{y \in \mathcal{V}, y \sim x}}w(x, y) F(x, t)\left[ G(y, t) - G(x, t)\right] \\
     &~~~+ \frac{1}{2}\sum_{y \in \mathcal{V}}  \sum_{{x \in \mathcal{V}, x \sim y}}w(x, y) F(y, t)\left[ G(x, t) - G(y, t)\right]\\
     &= - \frac{1}{2} \sum_{x, y\in \mathcal{V}}  w(x, y) \left[ F(y, t) - F(x, t)\right]\left[ G(y, t) - G(x, t)\right].
\end{align*}
\noindent \textbf{Corollary-1.}  If we choose $F(x, t) = 1$ for all $(x,t) \in \mathcal{V}\times [0, +\infty)$, then $$\sum\limits_{x \in \mathcal{V}}  \Delta G(x, t) = 0.$$
\noindent \textbf{Corollary-2.}  If we choose $F(x, t) =1/ {G(x, t)}$ for all $(x,t) \in \mathcal{V}\times [0, +\infty)$, then
\begin{align*}
    \sum_{x \in \mathcal{V}}\dfrac{1}{G(x, t)}  \Delta G(x, t) =  \frac{1}{2} \sum_{x, y \in \mathcal{V}} w(x, y) \dfrac{\left[ G(y, t) - G(x, t)\right]^2}{G(x, t) G(y, t) }.
\end{align*}
Now, we provide two theorems at disease-free equilibrium \textbf{D}$= (S_0, V_0, 0)$  and endemic equilibrium \textbf{E}$= (S_*, V_*, I_*)$ as following:

\begin{theorem}\label{thm}
\begin{enumerate}
    \item[(a)] The disease-free equilibrium $\textbf{D}$ of the model (\ref{model_reduce}) is globally asymptotically stable, when the basic reproduction number is less than unity i.e., $\mathcal{R}_0 < 1$
    \item[(b)] Consider $\eta = 0$. The endemic equilibrium $\textbf{E}$ of the model (\ref{model_reduce}) is globally asymptotically stable, when the basic reproduction number is greater than unity i.e., $\mathcal{R}_0 > 1$
\end{enumerate}
\end{theorem}
\noindent \textbf{Proof of Theorem \ref{thm}(a).} 
The disease-free equilibrium satisfies
\begin{equation}\label{disease_free_cond}
\left\lbrace
    \begin{aligned}
      \mu(1-r)N &=  (\mu+p)S_0  - \xi V_0 ,\\
    \mu r N   &=  -pS_0  + (\mu+\xi) V_0,
    \end{aligned}
    \right.
\end{equation}
Putting (\ref{disease_free_cond}) into (\ref{model_reduce}), we derive after some rearrangement in terms
\begin{equation}
\left\lbrace
    \begin{aligned}\label{model_dfe}
    \dfrac{dS}{dt} - \epsilon \Delta S &=  -(\mu+p)(S-S_0) - \beta (S-S_0) I - \beta S_0 I + \xi (V-V_0) + \eta N \\
    &~~~~- \eta(S-S_0) - \eta S_0 - \eta(V-V_0) - \eta V_0 - \eta I, \\
    \dfrac{dV}{dt} - \epsilon \Delta V &= p(S-S_0) - \beta \sigma (V-V_0) I - \beta \sigma V_0 - (\mu+\xi) (V-V_0),\\
    \dfrac{dI}{dt} - \epsilon \Delta I &= \beta (S-S_0) I + \beta S_0 I +  \beta \sigma (V-V_0) I + \beta \sigma V_0 I - (\mu+\gamma) I.
    \end{aligned}
    \right.
\end{equation}
To show the disease-free equilibrium is globally asymptotically stable, we define a Lyapunov function:

\begin{align}
    L(t) = \sum_{ x \in \mathcal{V}} \bigg[ \frac{1}{2S_0} (S-S_0)^2+ \frac{1}{2V_0} (V-V_0)^2+ I\bigg].
\end{align}
Therefore,  $L(t) \geq 0$ for all $t \geq 0$ and $L(t) = 0$ if and only if $(S, V, I) = (S_0, V_0, 0)$. Differentiating $L(t)$ with respect to `t', we get
\begin{align}\label{lya_der1}
    \frac{d}{dt} L(t) = \sum_{ x \in \mathcal{V}}\bigg[ \frac{S-S_0}{S_0} \frac{dS}{dt} + \frac{V-V_0}{V_0} \frac{dV}{dt}+ \frac{dI}{dt}\bigg].
\end{align}
Plugging all the derivatives from the system (\ref{model_dfe}) into (\ref{lya_der1}), we derive
\begin{align*}
    \frac{d}{dt} L(t) &= \sum_{ x \in \mathcal{V}}\bigg\lbrace \frac{S-S_0}{S_0} \Big[ -(\mu+p)(S-S_0) - \beta (S-S_0) I - \beta S_0 I + \xi (V-V_0) \\
    &~~~+ \eta N - \eta(S-S_0) - \eta S_0 - \eta(V-V_0) - \eta V_0 - \eta I + \epsilon \Delta S\Big] \\
    &~~~+ \frac{V-V_0}{V_0} \Big[p(S-S_0) - \beta \sigma (V-V_0) I - \beta \sigma V_0 - (\mu+\xi) (V-V_0)\\
    &~~~+\epsilon \Delta V  \Big]+ \Big[ \beta (S-S_0) I + \beta S_0 I +  \beta \sigma (V-V_0) I + \beta \sigma V_0 I - (\mu+\gamma) I + \epsilon \Delta I \Big] \bigg\rbrace.
\end{align*}
By using corollary-1 and corollary-2 of the Lemma-\ref{Green}, we obtain  $\sum_{ x \in \mathcal{V}} \epsilon \Delta I = 0$, 
\begin{align*}
    \sum_{ x \in \mathcal{V}} \frac{S-S_0}{S_0} \epsilon \Delta S &= \frac{\epsilon}{S_0}\sum_{ x \in \mathcal{V}} S \Delta S  - \epsilon\sum_{ x \in \mathcal{V}} \Delta S\\
    &= - \frac{\epsilon}{2S_0} \sum_{x, y \in \mathcal{V}} w(x, y) {\left[ S(y, t) - S(x, t)\right]^2} \leq 0, \quad \mbox{ for all } t \geq 0 \,\,\& \,\, \epsilon > 0.
\end{align*}
Similar way, one may derive
\begin{align*}
    \sum_{ x \in \mathcal{V}} \frac{V-V_0}{V_0} \epsilon \Delta V 
    &= - \frac{\epsilon}{2V_0} \sum_{x, y \in \mathcal{V}} w(x, y) {\left[ V(y, t) - V(x, t)\right]^2} \leq 0, \quad \mbox{ for all } t \geq 0 \,\,\& \,\, \epsilon > 0.
\end{align*}
Using the above results, we acquire
\begin{align*}
    \frac{d}{dt} L(t) &\leq  \sum_{ x \in \mathcal{V}}\bigg\lbrace \frac{S-S_0}{S_0} \Big[ -(\mu+p)(S-S_0) - \beta (S-S_0) I - \beta S_0 I + \xi (V-V_0) + \eta N\\
    &~~~- \eta(S-S_0) - \eta S_0 - \eta(V-V_0) - \eta V_0 - \eta I \Big] + \frac{V-V_0}{V_0} \Big[p(S-S_0)\\
    &~~~  - \beta \sigma (V-V_0) I - \beta \sigma V_0 I - (\mu+\xi) (V-V_0)\Big]\\
    &~~~+ \Big[ \beta (S-S_0) I  +  \beta \sigma (V-V_0) I\Big] +  \Big[\beta (S_0+  \sigma V_0) - (\mu+\gamma) \Big]I\bigg\rbrace.
\end{align*}
If $\mathcal{R}_0 < 1$, implies $\dfrac{\beta(S_0+\sigma V_0)}{\gamma +\mu } < 1$ i.e., $\beta(S_0+\sigma V_0) < (\gamma+\mu)$. Hence, the last term is a negative number. The above relation can be written as
\begin{equation}
\begin{aligned}\label{lap_inq}
    \frac{d}{dt} L(t) &\leq  \sum_{ x \in \mathcal{V}} \bigg\lbrace T(S, V)  +   \Big[  - \beta \frac{(S-S_0)^2}{S_0} I    - \eta \frac{(S-S_0)^2}{S_0} - \beta \sigma \frac{(V-V_0)^2}{V_0}I\Big] \\
    &~~~-\eta {(S-S_0)}\Big[ 1+\frac{V_0}{S_0} -\frac{N}{S_0}  \Big]\bigg\rbrace\\
    &=  \sum_{ x \in \mathcal{V}} \bigg\lbrace T(S, V)  +   \Big[  - \beta \frac{(S-S_0)^2}{S_0} I    - \eta \frac{(S-S_0)^2}{S_0} - \beta \sigma \frac{(V-V_0)^2}{V_0}I\Big]\bigg\rbrace,
\end{aligned}
\end{equation}
 since $S_0 +V_0 = N$ and  where $T(S, V) = -\Big[\frac{\mu+p}{S_0} + \frac{\eta}{S_0}\Big](S-S_0)^2 + \Big[\frac{\xi}{S_0} - \frac{\eta}{S_0} + \frac{p}{V_0}\Big](S-S_0)(V-V_0)  -\frac{\mu+\xi}{V_0}(V-V_0)^2  $. We rewrite $T(S, V)$ as  $T(S, V) = A(S-S_0)^2 + B(S-S_0)(V-V_0) + C(V-V_0)^2 $, where
\begin{align}
    A =-\Big[\frac{\mu+p}{S_0} + \frac{\eta}{S_0}\Big], \quad B =  \frac{\xi}{S_0} - \frac{\eta}{S_0} + \frac{p}{V_0}, \quad C = -\frac{\mu+\xi}{V_0}. 
\end{align}
Now, we show that the function $T(S, V)$ is negative definite with respect to $S=S_0$ and $V = V_0$. To show that, we need to verify $B^2 - 4AC < 0$. We calculate
\begin{align*}
    B^2 - 4AC &= \Big[\frac{\xi}{S_0} - \frac{\eta}{S_0} + \frac{p}{V_0}\Big]^2 - 4 \Big[\frac{\mu+p}{S_0} + \frac{\eta}{S_0}\Big] \Big[ \frac{\mu+\xi}{V_0} \Big]\\
    &= \frac{1}{S_0^2} \Big\lbrace\Big[{\xi} - {\eta} + \frac{p S_0}{V_0}\Big]^2 - 4 ({\mu+p} + {\eta})(\mu+\xi) \frac{S_0}{V_0}\Big\rbrace\\
    &\leq  \frac{1}{S_0^2} \Big\lbrace\Big[{\xi}  + \frac{p S_0}{V_0}\Big]^2 - 4 ({\mu+p})(\mu+\xi) \frac{S_0}{V_0}\Big\rbrace  - 4  {\eta}(\mu+\xi) \frac{1}{S_0V_0}.
\end{align*}
Compute ${S_0}/{V_0}  = (\xi+(1-r)\mu)/(p+r\mu)$. We have $\frac{p}{p+r\mu} < 1$ and $\frac{p+\mu}{p+r\mu} > 1$, since $0 < r < 1$. Using these results, the above becomes
\begin{align*}
    B^2 - 4AC &\leq  \frac{1}{S_0^2} \Big\lbrace (2\xi +(1-r)\mu)^2 - 4 (\mu+\xi) (\xi +(1-r)\mu)\Big\rbrace  - 4  {\eta}(\mu+\xi) \frac{1}{S_0V_0}.
\end{align*}
With some straight calculations and cancellation of terms, we derive
\begin{align*}
    B^2 - 4AC &\leq  \frac{1}{S_0^2} \Big\lbrace \mu^2(1-r)[(1-r)-4] - 4 \xi \mu \Big\rbrace  - 4  {\eta}(\mu+\xi) \frac{1}{S_0V_0}.
\end{align*}
Again since, $0 < 1-r <1 $, $(1-r) < 4$. Thus, we obtain $B^2 - 4AC < 0$. Finally, from (\ref{lap_inq}) it is obvious that 
\begin{align*}
    \frac{d}{dt} L(t) &\leq    \sum_{ x \in \mathcal{V}} \bigg\lbrace T(S, V)  +   \Big[  - \beta \frac{(S-S_0)^2}{S_0} I    - \eta \frac{(S-S_0)^2}{S_0} - \beta \sigma \frac{(V-V_0)^2}{V_0}I\Big]\bigg\rbrace < 0, 
\end{align*}
$\mbox{ for all } S, V, I \geq 0.$ By applying the LaSalle's Lyapunov Principle, see \cite{Lyapunov, Lasalle}, 
\begin{align*}
    \lim_{ t \to \infty} (S, V, I) = (S_0, V_0, 0), \quad \mbox{ uniformly in } x \in \mathcal{V}.
\end{align*}
Hence, disease-free equilibrium, \textbf{D}, is globally asymptotically stable provided $\mathcal{R}_0 < 1$.\\[2ex]
\noindent \textbf{Proof of Theorem \ref{thm}(b).} 
The endemic equilibrium point satisfies
\begin{equation}
\left\lbrace
    \begin{aligned}\label{endemic_cond}
    \mu+p &= \mu(1-r)\frac{N}{S_*}  - \beta  I_* + \xi \frac{V_*}{S_*},\\
    \mu+\xi &=\mu r \frac{N}{V_*}  + p \frac{S_*}{V_*} - \beta \sigma I_*, \\
    \mu+\gamma &= \beta S_* +  \beta \sigma V_*.
    \end{aligned}
    \right.
\end{equation}
For showing the endemic equilibrium \textbf{E} is globally asymptotically stable, we consider the following Lyapunov function:
\begin{align}
L_*(t) = \sum_{x \in \mathcal{V}} \Big[ \Big(S - S_* - S_* \ln \dfrac{S}{S_*}\Big) + \Big(V - V_* - V_* \ln \dfrac{V}{V_*}\Big) + \Big(I - I_* - I_* \ln \dfrac{I}{I_*}\Big) \Big]
\end{align}
Compute its derivative with respect to 't' and then apply \eqref{model_reduce} in the case where $\eta = 0$, 
\begin{align}\label{Lder}\notag
    \dfrac{d}{dt} L_*(t) &= \sum_{x \in \mathcal{V}} \Big[\Big(1- \dfrac{S_*}{S}\Big)\Big(\mu(1-r)N -(\mu+p)S - \beta S I + \xi V + \epsilon \Delta S\Big)\\
    &~~~+\Big(1- \dfrac{V_*}{V}\Big)\Big(\mu r N + pS - \beta \sigma V I - (\mu+\xi) V + \epsilon \Delta V\Big) \\
    &~~~+ \Big(1- \dfrac{I_*}{I}\Big)\big( \beta S I +  \beta \sigma V I - (\mu+\gamma) I  + \epsilon \Delta I\Big)\Big].\notag
\end{align}
Substituting (\ref{endemic_cond}) into (\ref{Lder}), we derive
\begin{align}\label{Lder1}\notag
    \dfrac{d}{dt} L_*(t) &=  \sum_{x \in \mathcal{V}} \Big[ \mu(1-r)N \Big(2 - \dfrac{S_*}{S} - \dfrac{S}{S_*}\Big) + \mu r N \Big(2 - \dfrac{V_*}{V} - \dfrac{V}{V_*}\Big)\\
    &~~~  + \epsilon  \Big(1- \dfrac{S_*}{S}\Big) \Delta S  + \epsilon  \Big(1- \dfrac{V_*}{V}\Big) \Delta V + \epsilon  \Big(1- \dfrac{I_*}{I}\Big) \Delta I\\
    &~~~ + \xi V_*\Big(\dfrac{V}{V_*} - \dfrac{S}{S_*} \Big)\Big( 1 - \dfrac{S_*}{S}\Big) + p S_*\Big(\dfrac{S}{S_*} - \dfrac{V}{V_*} \Big)\Big( 1 - \dfrac{V_*}{V} \Big)  \Big].\notag
\end{align}
By using corollary-1 and corollary-2 of the Lemma-\ref{Green}, we obtain  
\begin{align*}
    \sum_{ x \in \mathcal{V}} \epsilon  \Big(1- \dfrac{S_*}{S}\Big) \Delta S &= \epsilon\sum_{ x \in \mathcal{V}} \Delta S -  {\epsilon S_* }\sum_{ x \in \mathcal{V}} \dfrac{1}{S}\Delta S  \\
    &= - \frac{\epsilon S_*}{2} \sum_{x, y \in \mathcal{V}} w(x, y) \dfrac{\left[ S(y, t) - S(x, t)\right]^2}{S(x, t) S(y, t) } \leq 0,
\end{align*}
$\mbox{for all }  t \geq 0, \,\,\mbox{ and }\,\, \epsilon > 0$. Similar way, one may derive
\begin{align*}
     \sum_{ x \in \mathcal{V}} \epsilon  \Big(1- \dfrac{V_*}{V}\Big) \Delta V
    &= - \frac{\epsilon V_*}{2} \sum_{x, y \in \mathcal{V}} w(x, y) \dfrac{\left[ V(y, t) - V(x, t)\right]^2}{V(x, t) V(y, t) } \leq 0,\\
   \sum_{ x \in \mathcal{V}} \epsilon  \Big(1- \dfrac{I_*}{I}\Big) \Delta I
    &= - \frac{\epsilon I_*}{2} \sum_{x, y \in \mathcal{V}} w(x, y) \dfrac{\left[ I(y, t) - I(x, t)\right]^2}{I(x, t) I(y, t) } \leq 0, 
\end{align*}
$\mbox{for all }t \geq 0, \,\,\mbox{ and }\,\, \epsilon > 0$. Using the above three results into \eqref{Lder1}, we obtain
\begin{align}\label{g_fun}
    \dfrac{d}{dt} L_*(t) &\leq  \sum_{x \in \mathcal{V}} g(S, V), 
\end{align}
where 
\begin{align*}
    g(S, V) &=  \mu(1-r)N \Big(2 - \dfrac{S_*}{S} - \dfrac{S}{S_*}\Big) + \mu r N \Big(2 - \dfrac{V_*}{V} - \dfrac{V}{V_*}\Big)\\
    &~~~ + \xi V_*\Big(\dfrac{V}{V_*} - \dfrac{V S_*}{V_*S} - \dfrac{S_*}{S} + 1 \Big) + p S_*\Big(\dfrac{S}{S_*} - \dfrac{S V_*}{S_*V} - \dfrac{V}{V_*} + 1 \Big).
\end{align*}
Now, to show the function $g$ is negative,  we consider three cases: $pS_* > \xi V_*,\, pS_* = \xi V_*$ and $pS_* < \xi V_*$:
\begin{enumerate}
    \item \textbf{When $pS_* = \xi V_*$:}
\begin{align*}
    g(S, V) &=  \mu(1-r)N \Big(2 - \dfrac{S_*}{S} - \dfrac{S}{S_*}\Big) + \mu r N \Big(2 - \dfrac{V_*}{V} - \dfrac{V}{V_*}\Big)\\
    &~~~ + \xi V_*\Big(2 - \dfrac{V S_*}{V_*S} - \dfrac{V_* S}{V S_*} \Big).
\end{align*}
 \item \textbf{When $pS_* > \xi V_*$:} from the first equation of \eqref{endemic_cond}, we have
 $\mu(1- r) N = (\mu + \beta I_*)S_* + pS_* - \xi V_*$. Substituting this into $g$, we obtain
\begin{align*}
    g(S, V) &=  \mu r N \Big(2 - \dfrac{V_*}{V} - \dfrac{V}{V_*}\Big) + (\mu + \beta I_*)S_*\Big(2 - \dfrac{S_*}{S} - \dfrac{S}{S_*}\Big)\\
    &~~~ + (pS_* - \xi V_*) \Big(3 - \dfrac{V}{V_*} - \dfrac{S_*}{S} - \dfrac{V_* S}{V S_*} \Big) + \xi V_*\Big( 2 - \dfrac{V_* S}{V S_*} - \dfrac{V S_*}{V_*S} \Big).
\end{align*}
 \item \textbf{When $pS_* < \xi V_*$:} from the second equation of \eqref{endemic_cond}, we have
$\mu r N = (\mu + \sigma \beta I_*)V_* - pS_* + \xi V_*$. Substituting this into $g$, we obtain
\begin{align*}
    g(S, V) &=  \mu(1-r)N \Big(2 - \dfrac{S_*}{S} -\dfrac{S}{S_*}\Big) + (\mu + \sigma \beta I_*)V_*\Big(2 - \dfrac{V_*}{V} - \dfrac{V}{V_*}\Big)\\
    &~~~ + (\xi V_* - pS_*) \Big(3 - \dfrac{S}{S_*} - \dfrac{V_*}{V} - \dfrac{V S_*}{V _*S} \Big) + p S_*\Big( 2 - \dfrac{V_* S}{V S_*} - \dfrac{V S_*}{V_*S} \Big).
\end{align*}
\end{enumerate}
\mbox{}\noindent Applying arithmetic mean $\geq$ geometric mean, it is clear that $g$ is a negative definite with respect to $S=  S_*$ and $V = V_*$.
Using the all the above results into (\ref{g_fun}), we have $ \frac{d}{dt} L_*(t) \leq 0$ for all $t \geq 0$.  By applying the LaSalle's Lyapunov Principle, see \cite{Lyapunov, Lasalle}, 
\begin{align*}
    \lim_{ t \to \infty} (S, V, I) = (S_*, V_*, I_*), \quad \mbox{ uniformly in } x \in \mathcal{V}.
\end{align*}
Hence, endemic equilibrium, \textbf{E}, is globally asymptotically stable provided $\mathcal{R}_0 > 1$.\\
\mbox{}
\begin{figure}[h!]
    \centering
\mbox{}\hspace*{-15pt}    \includegraphics[width=0.8\textwidth]{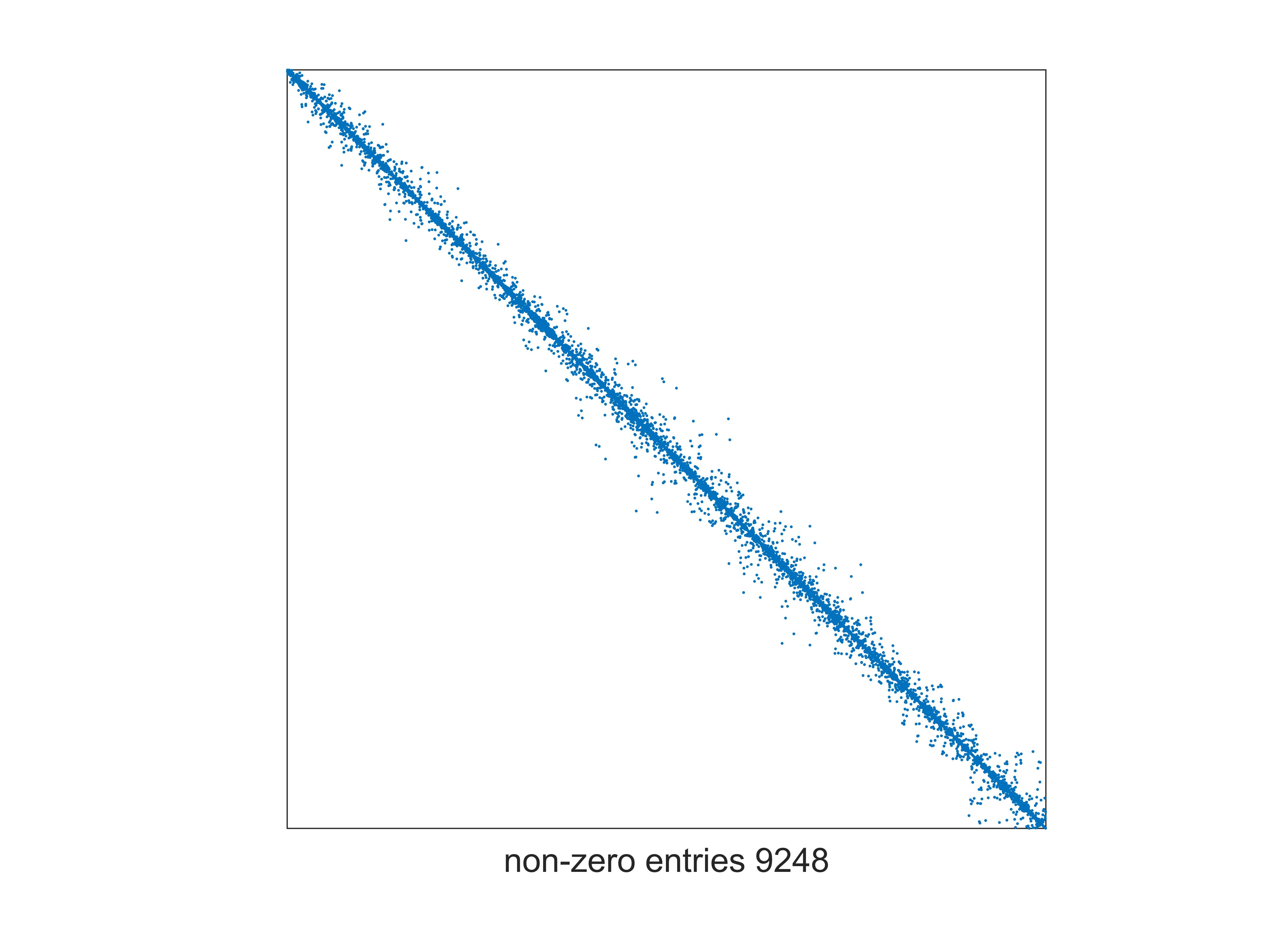}
    \caption{Sparsity of the Laplacian matrix corresponding the Minnesota road network. The   graphical structure of the network is given in Fig. \ref{fig:graph}.}
    \label{fig:sparsity}
\end{figure}
\section{Numerical Results and Discussions}\label{numerical}
\begin{figure}[h!]
    \centering
\hspace*{-10pt}\begin{tabular}{c c}
      {\includegraphics[width=80mm]{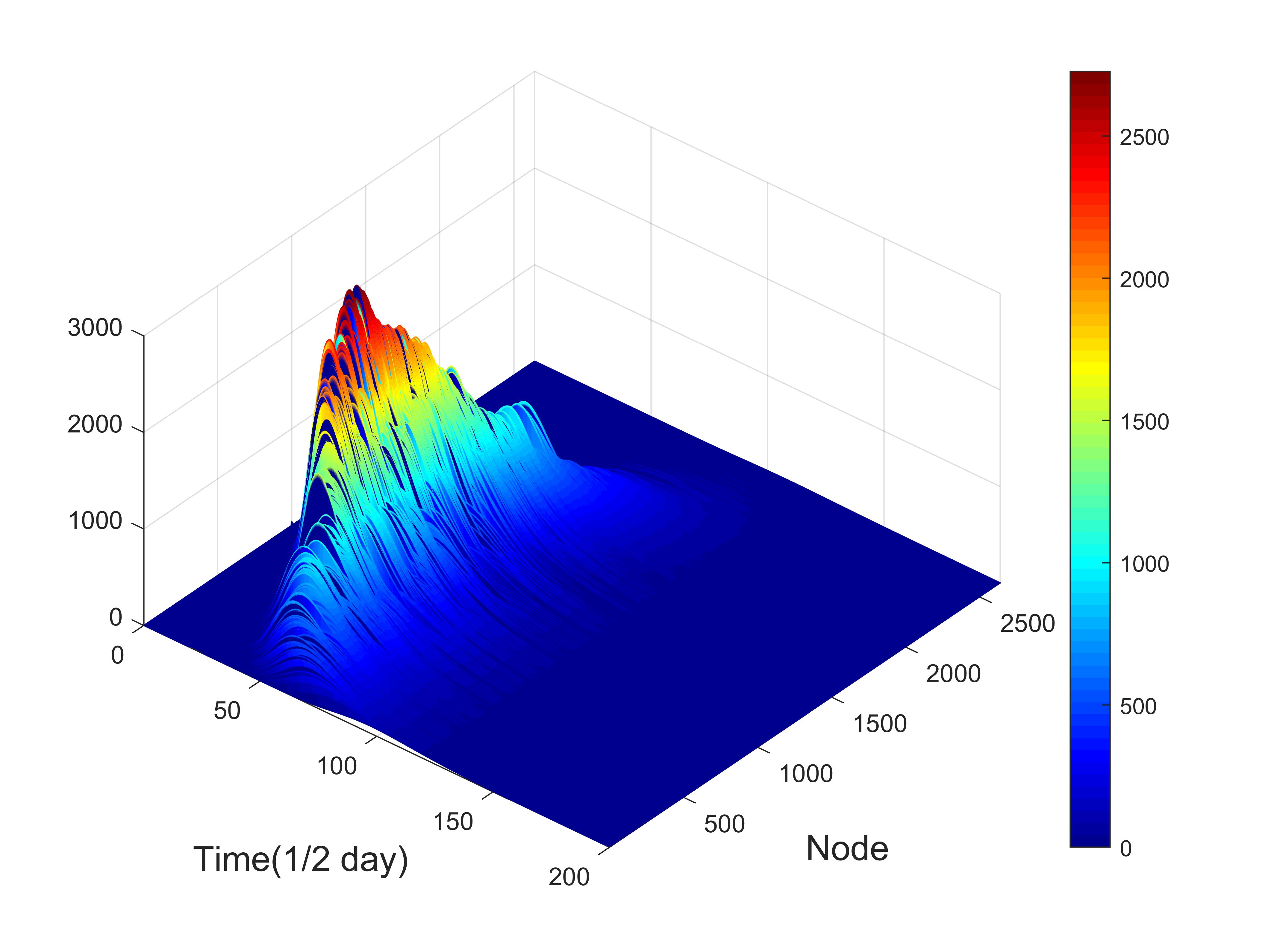}} &
      {\includegraphics[width=80mm]{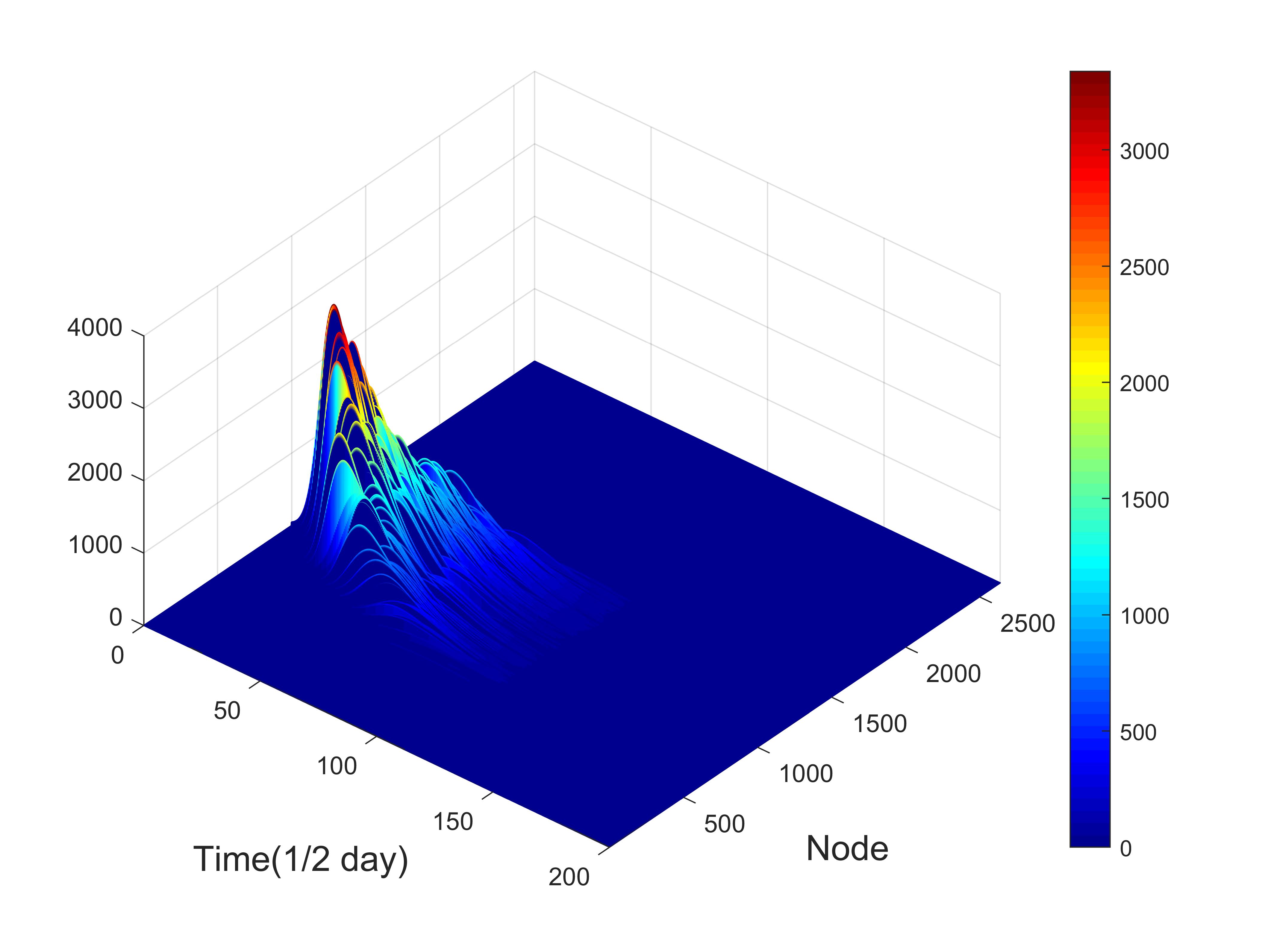}} \\
      \small Solution of $I(x, t)$ when $\epsilon = 1$ &  Solution of $I(x, t)$ when $\epsilon = 0.1$ \\
      {\includegraphics[width=80mm]{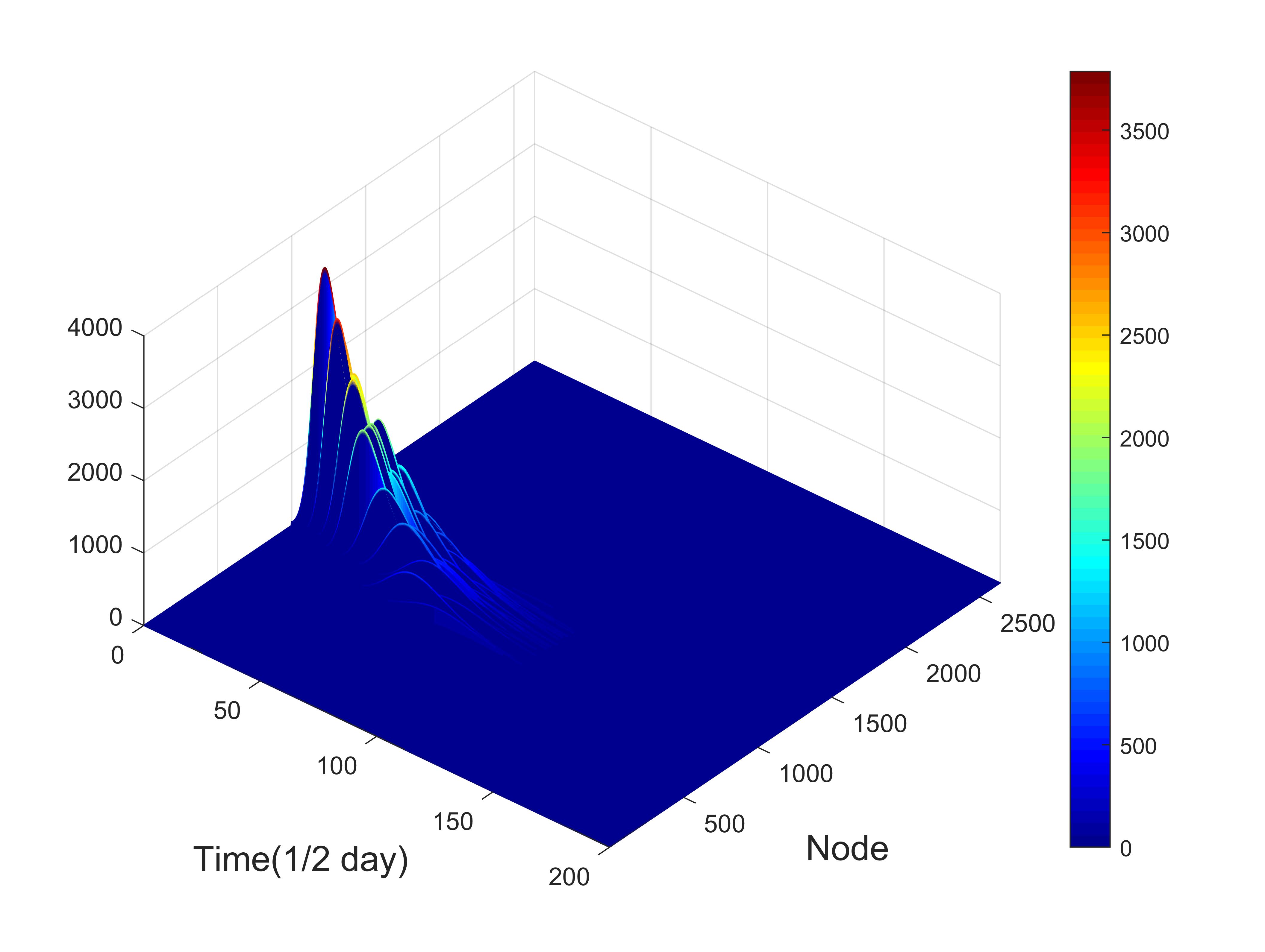}} &
      {\includegraphics[width=80mm]{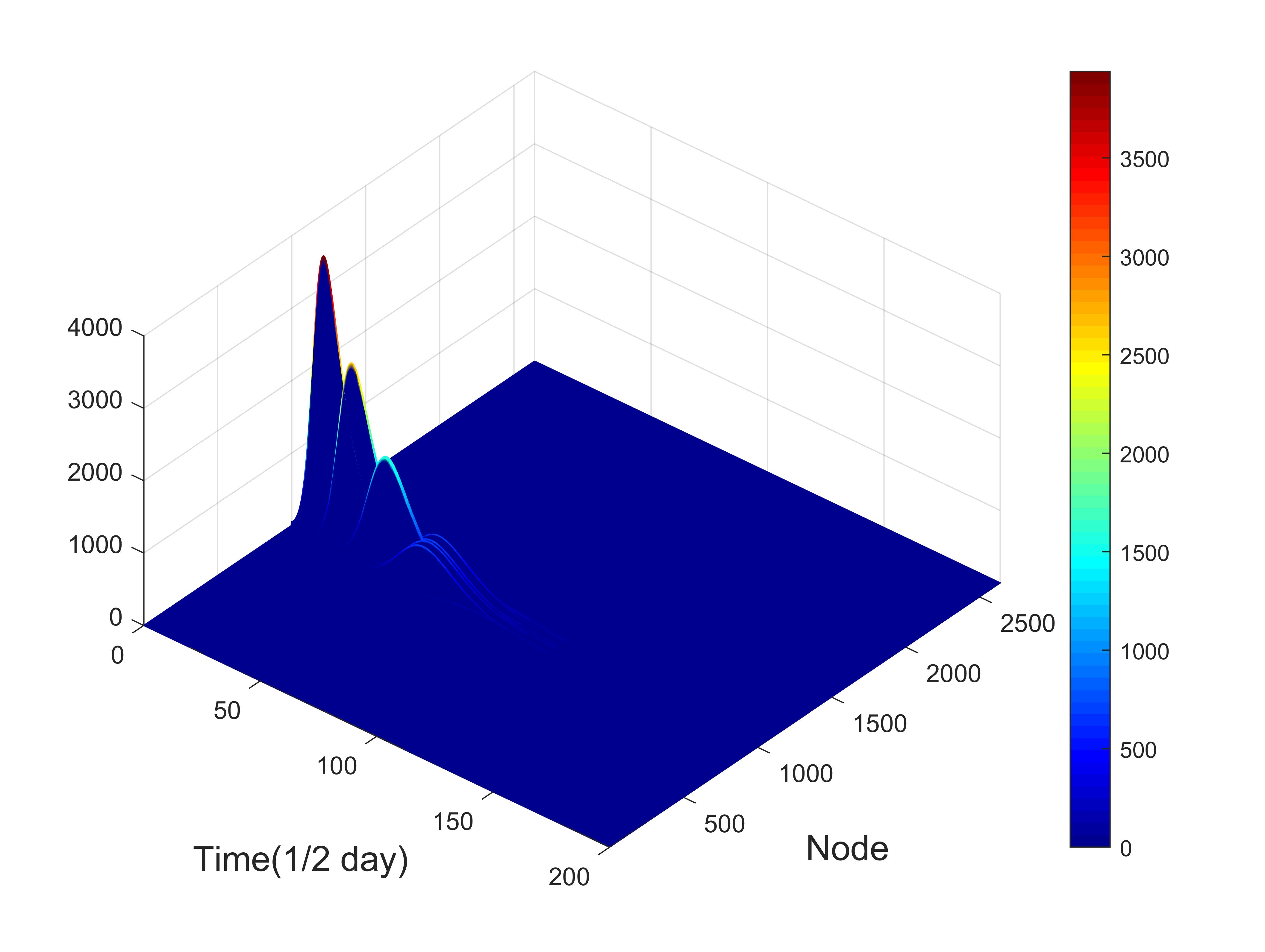}} \\
      \small Solution of $I(x, t)$ when $\epsilon = 0.01$ & Solution of $I(x, t)$ when $\epsilon = 0.001$      
\end{tabular}
    \caption{The figure portrays that the solution of $I(x, t)$ while vaccination rate ($p$) is fixed at 0.02 and a fraction of newborn vaccinated is fixed at $r=0.001$. Moreover, $\xi = \eta = \sigma = 0.0$. }
    \label{fig:vac_vs_mig}
\end{figure}

\noindent This section is dedicated to all numerical results and their discussions of the model (\ref{model}). The graph $\mathcal{V}$ is collected from UF Sparse Matrix Collection, see \cite{network}. It is a Minnesota road network located in the upper mid-western region of the United States. It consists of 2642 nodes and the weight between any two adjacent nodes is unity., i.e., $w(x, y) = 1$. The size of the Laplacian matrix corresponding to the graph is $2642\times2642$, which is a sparse matrix with 9248 non-zero entries. The sparsity of the matrix is illustrated in Figure \ref{fig:sparsity}, where all the non-zero entries have occurred near the diagonal. This means there is no direct road between two nodes that are placed sufficiently far from each other. To make numerical calculations simple, we assumed each node has the same population of 10,000 and natural death rate $\mu = 0.001$. These two parameters are fixed throughout the paper. Other parameters will be mentioned based on various analyses when needed in the following:
\begin{figure}[hbt!]
    \centering
\begin{tabular}{c }
      {\includegraphics[width=12cm, height=10cm]{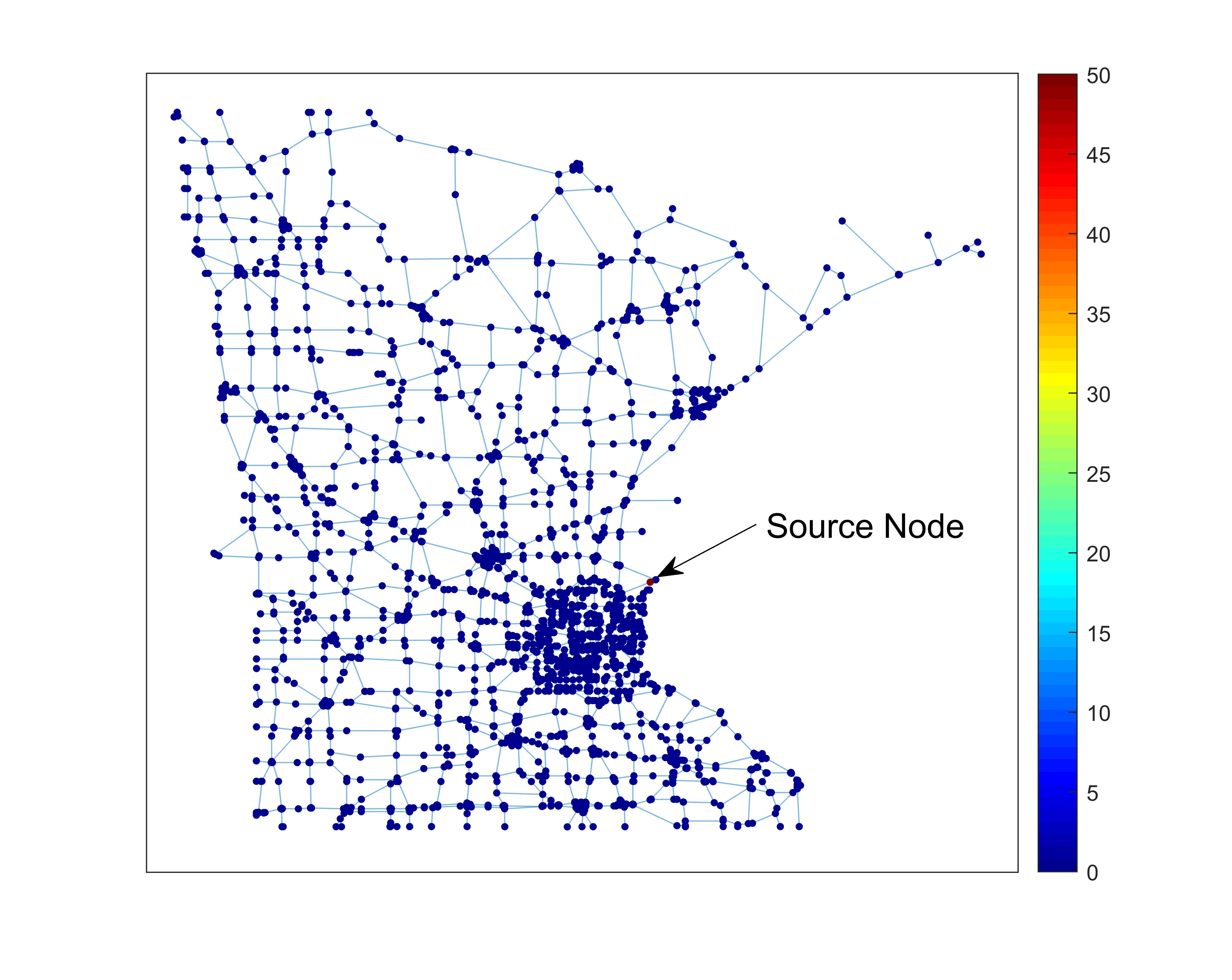}} \\
      \small Evolution of $I(x, t)$ at time $t$=$0$ and node $x$=$990$ \\
      {\includegraphics[width=12cm, height=10cm]{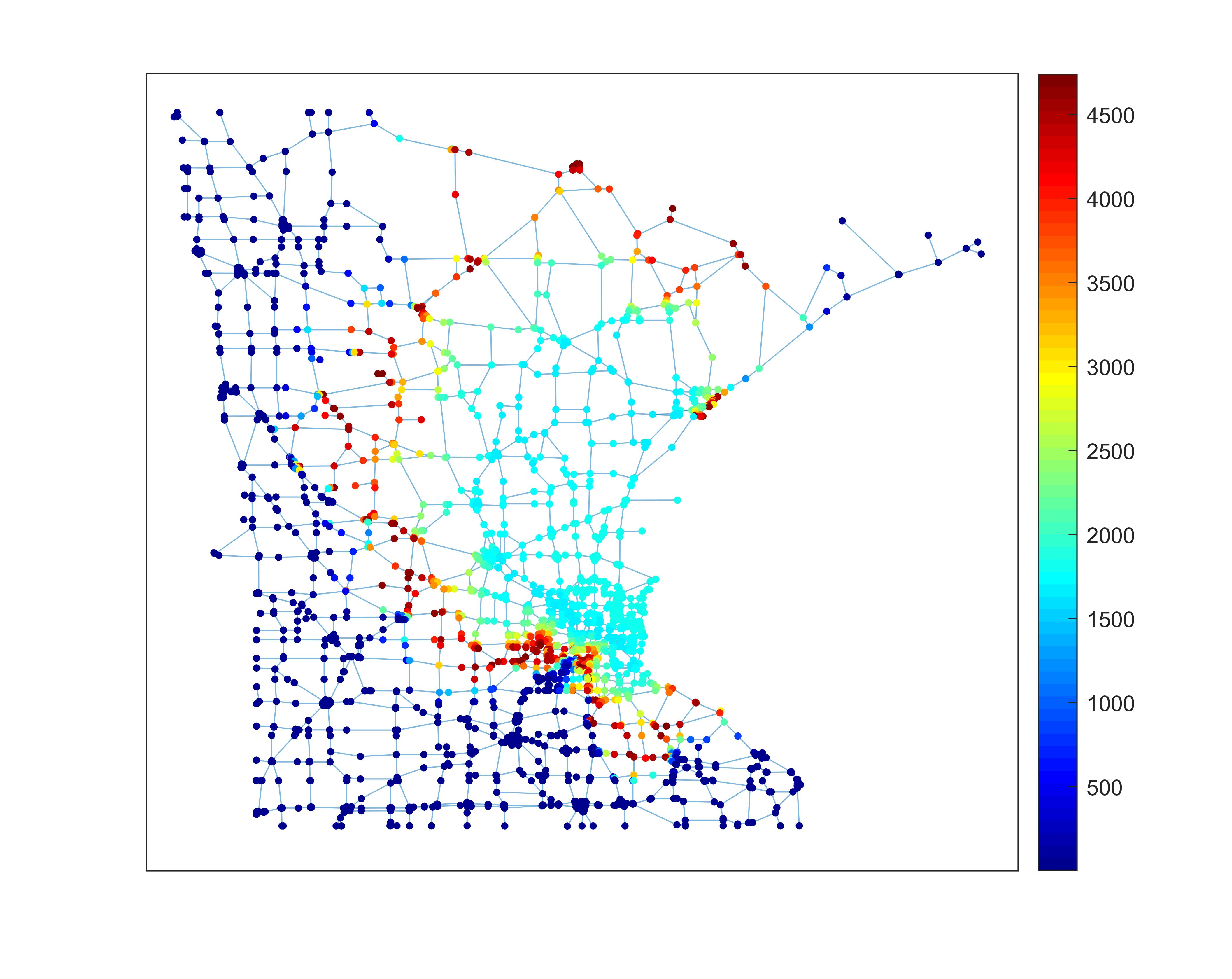}} \\
       Evolution of $I(x, t)$ at time $t$=$70$ days and node $x$=$990$ \\
\end{tabular}
    \caption{The figure depicts the solution of infected individual $I(x, t)$, when the parameters are as follows: $\beta = 0.65/N, r =0.01, p = 0.002, \xi = 0.02, \eta = 0.03, \gamma = 0.11$. Furthermore, we assumed that the migration parameter rate or population mobility rate $\epsilon = 0.05$. }
    \label{fig:graph}
\end{figure}
we consider initially at $t = 0$ the infection starts from node 990, which we are calling the source node,  with 50 infected individuals, 10 vaccinated individuals, and all other nodes have zero infected and vaccinated individuals in the network. The source node 990 is pointed out in Figure \ref{fig:graph} (above). The evolution of the solution $I(x,t)$ is presented in Figure \ref{fig:graph} (below)  at $t = 70$ days, which displays a spatially non-homogeneous behavior. The parameters are given corresponding to this result in the caption of Figure \ref{fig:graph}. Here we have chosen population mobility rate $\epsilon = 0.05$, which restricts 95\% of the population movement. Moreover, it is observed that as $\epsilon$ tends to one the peak of infection occurs at each node simultaneously, which results in an epidemic in the whole network.\\[2ex]
\noindent In Figure \ref{fig:vac_vs_mig}, we discussed the effect of population mobility on limited vaccination. Control the spread of disease is more cost-effective by increasing vaccination. If the vaccination is limited then we need some strategies to control infections. For this, we set $\xi = \sigma = 0$, i.e., the vaccine is completely effective, and $\eta = 0$, i.e., no loss of immunity who are in the removal compartment. The vaccination rate and fraction of newborns vaccinated are fixed at $p = 0.02$ and $r = 0.001$, respectively. For $\epsilon = 1$, there is no restriction on movement from one node to another, the disease starts spreading from the source node to another, and because of vaccination, the spread decays as the distance from the source node increases as shown in Figure \ref{fig:vac_vs_mig}(upper-left). From the sub-figures(upper-right, lower-left and lower-right), when $\epsilon = 0.1, 0.01, 0.001$, as we increase the restriction of population mobility on fixed vaccination rate, the infection decays drastically as the distance increases from the source node.\\[2ex]
\begin{figure}[!ht]
    \centering
\hspace*{-10pt}\begin{tabular}{c c}
      {\includegraphics[width=80mm]{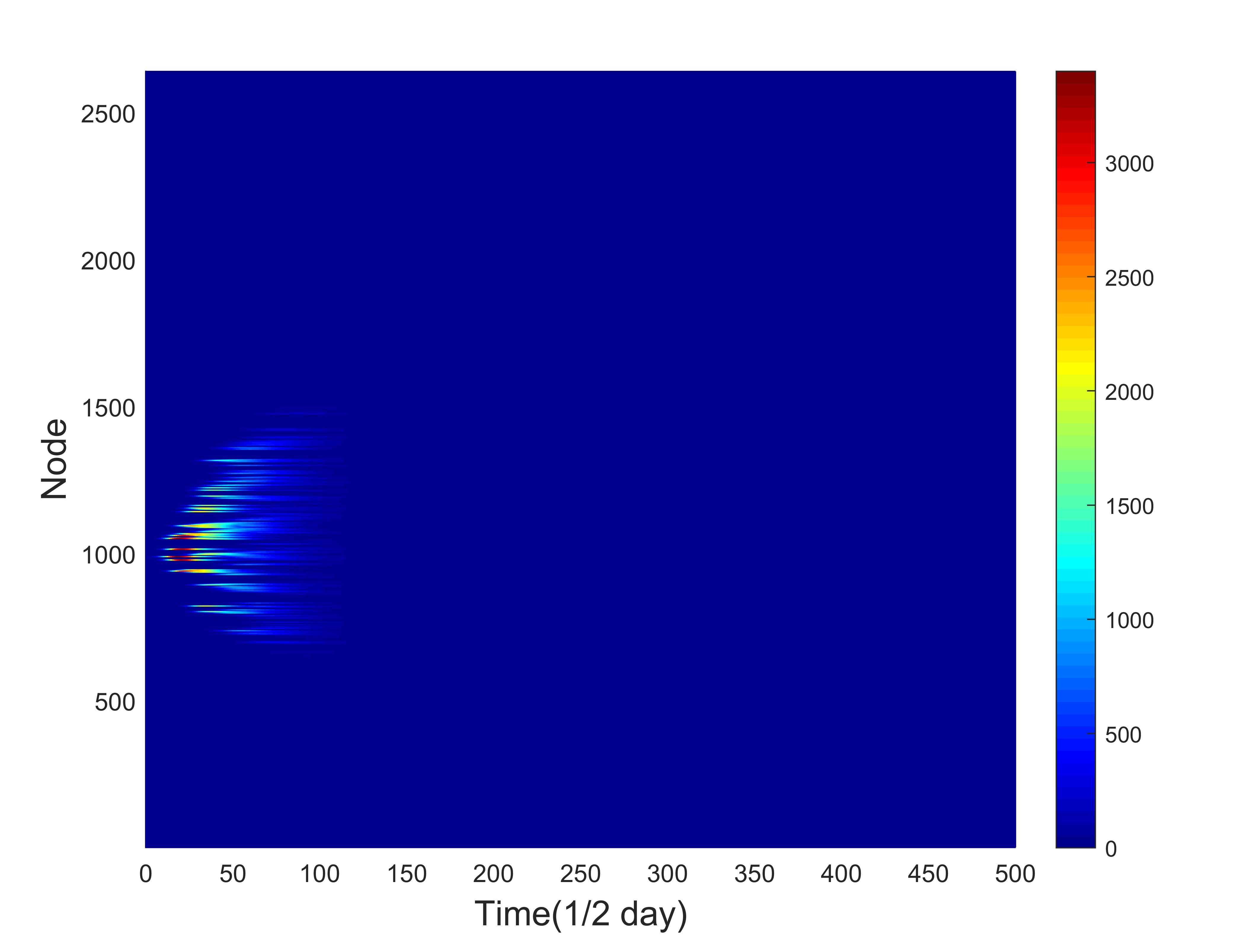}} &
      {\includegraphics[width=80mm]{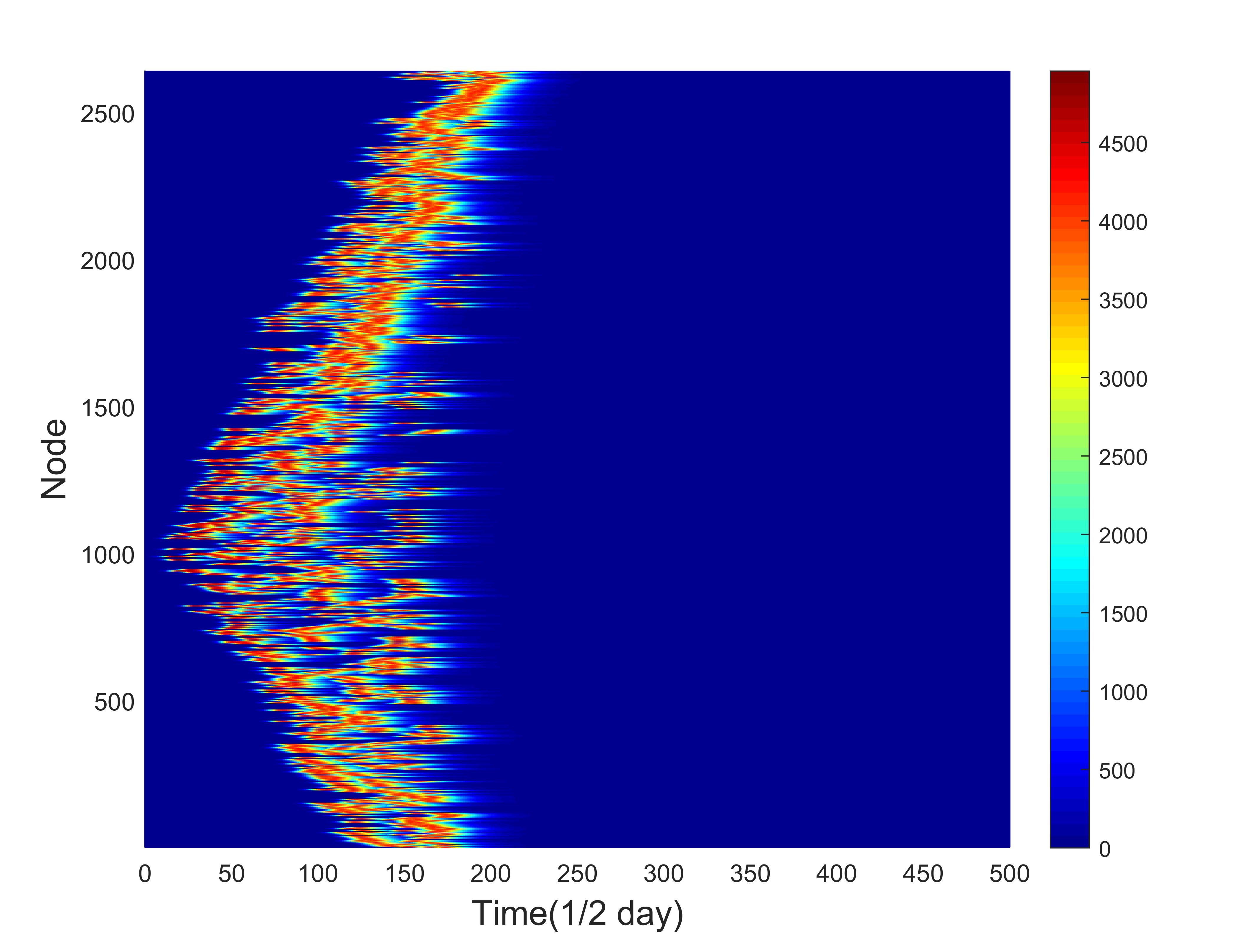}} \\
      \small Solution of $I(x, t)$ when $\sigma = \eta = \xi = 0$ &  Solution of $I(x, t)$ when $\sigma = 0.7$ \\
      {\includegraphics[width=80mm]{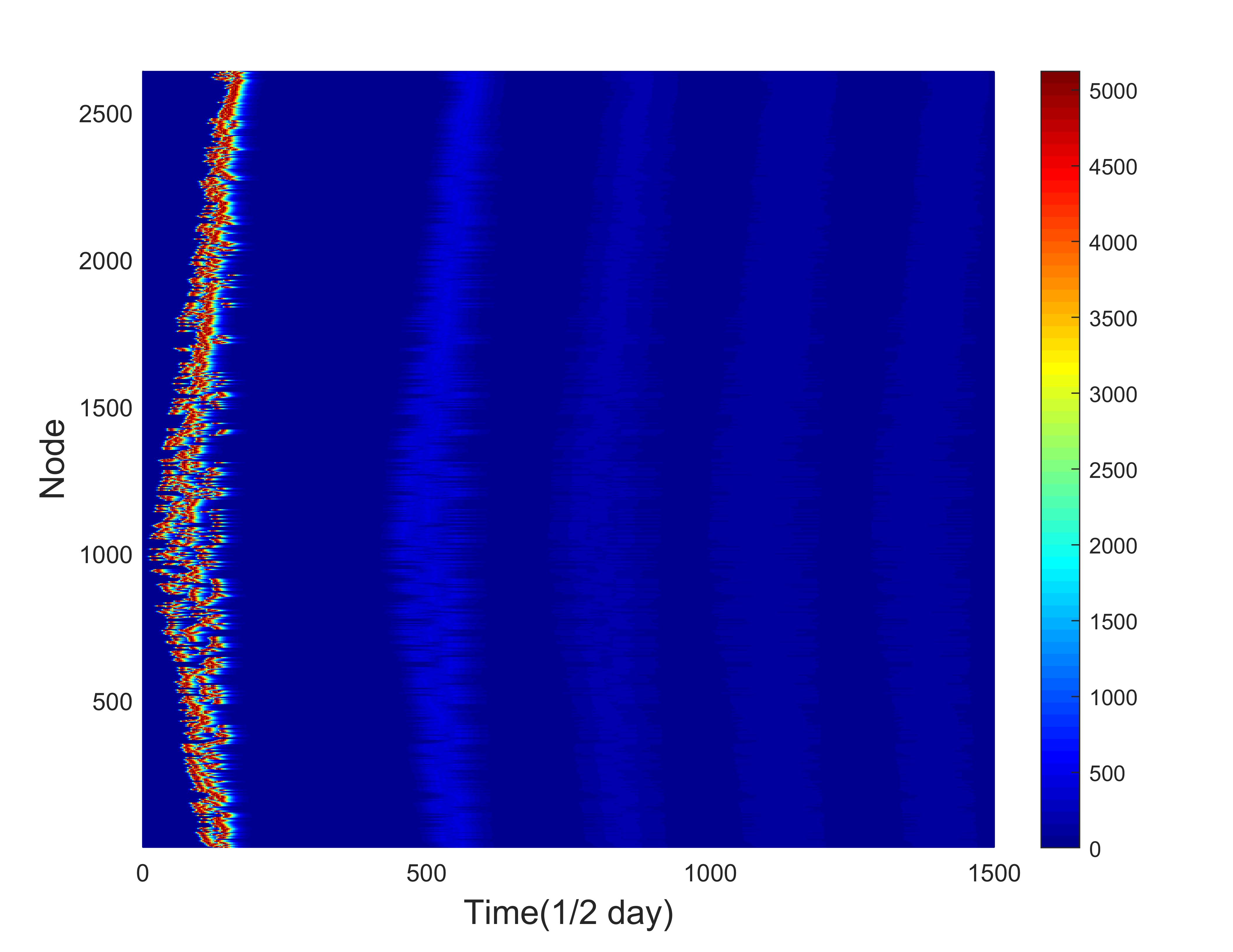}} &
      {\includegraphics[width=80mm]{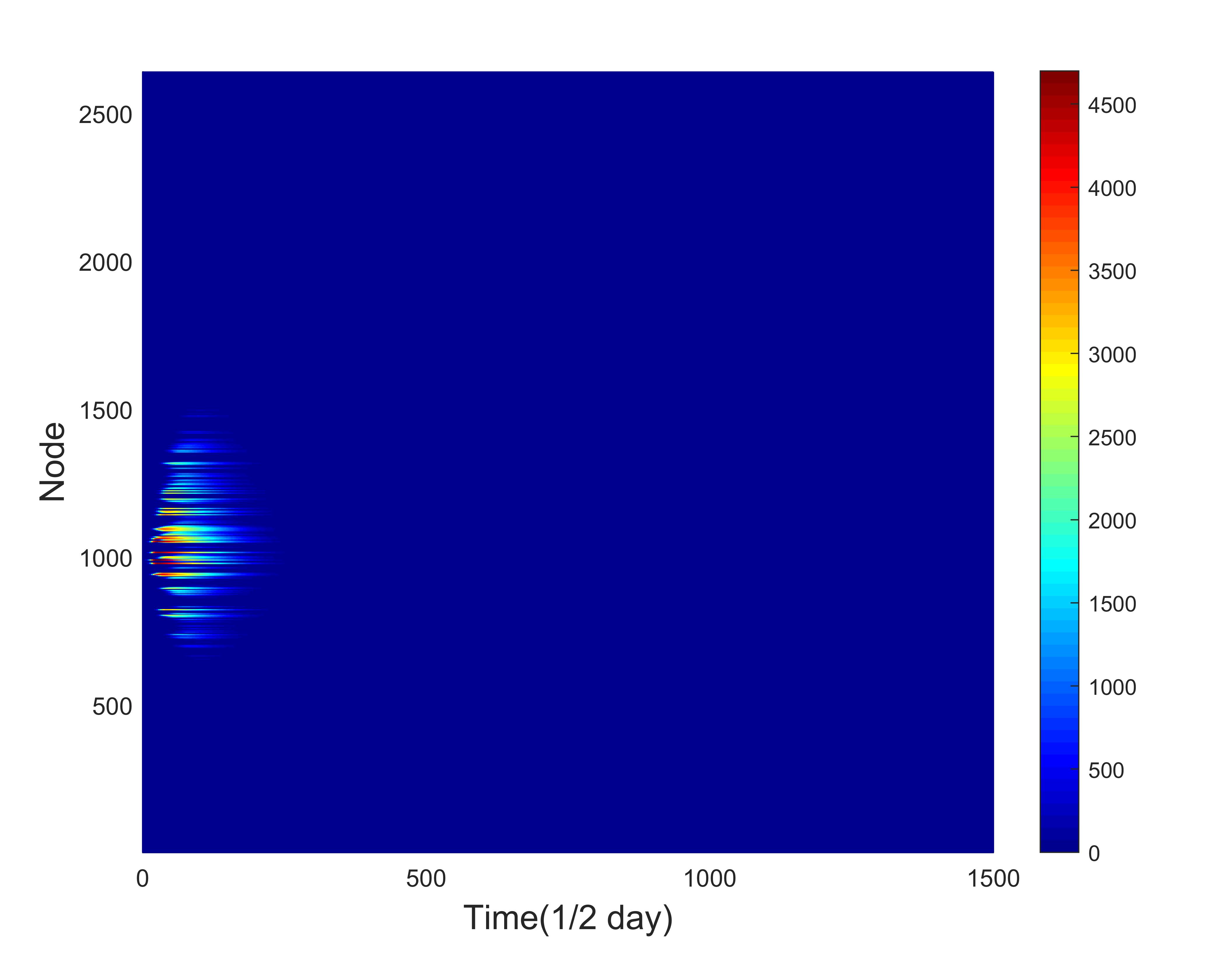}} \\
      \small Solution of $I(x, t)$ when $\xi = 0.7$ & Solution of $I(x, t)$ when $\eta = 0.7$      
\end{tabular}
    \caption{The figure portraits that the solution of $I(x, t)$ while the population mobility rate is fixed at $\epsilon = 0.05$.}
    \label{fig:var}
\end{figure}
\begin{figure}[h!]
    \centering
\hspace*{-10pt}\begin{tabular}{c c}
      {\includegraphics[width=80mm]{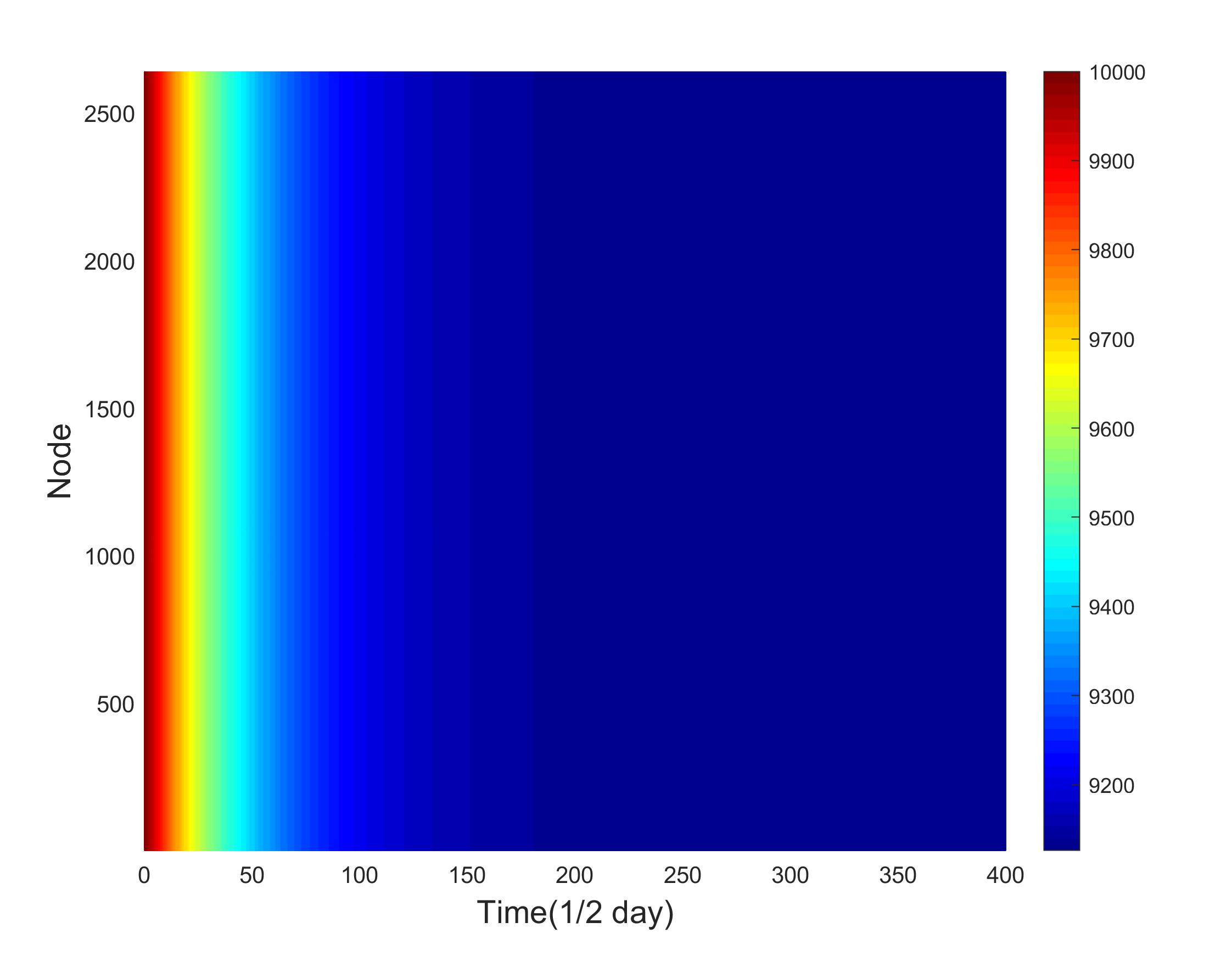}} &
      {\includegraphics[width=80mm]{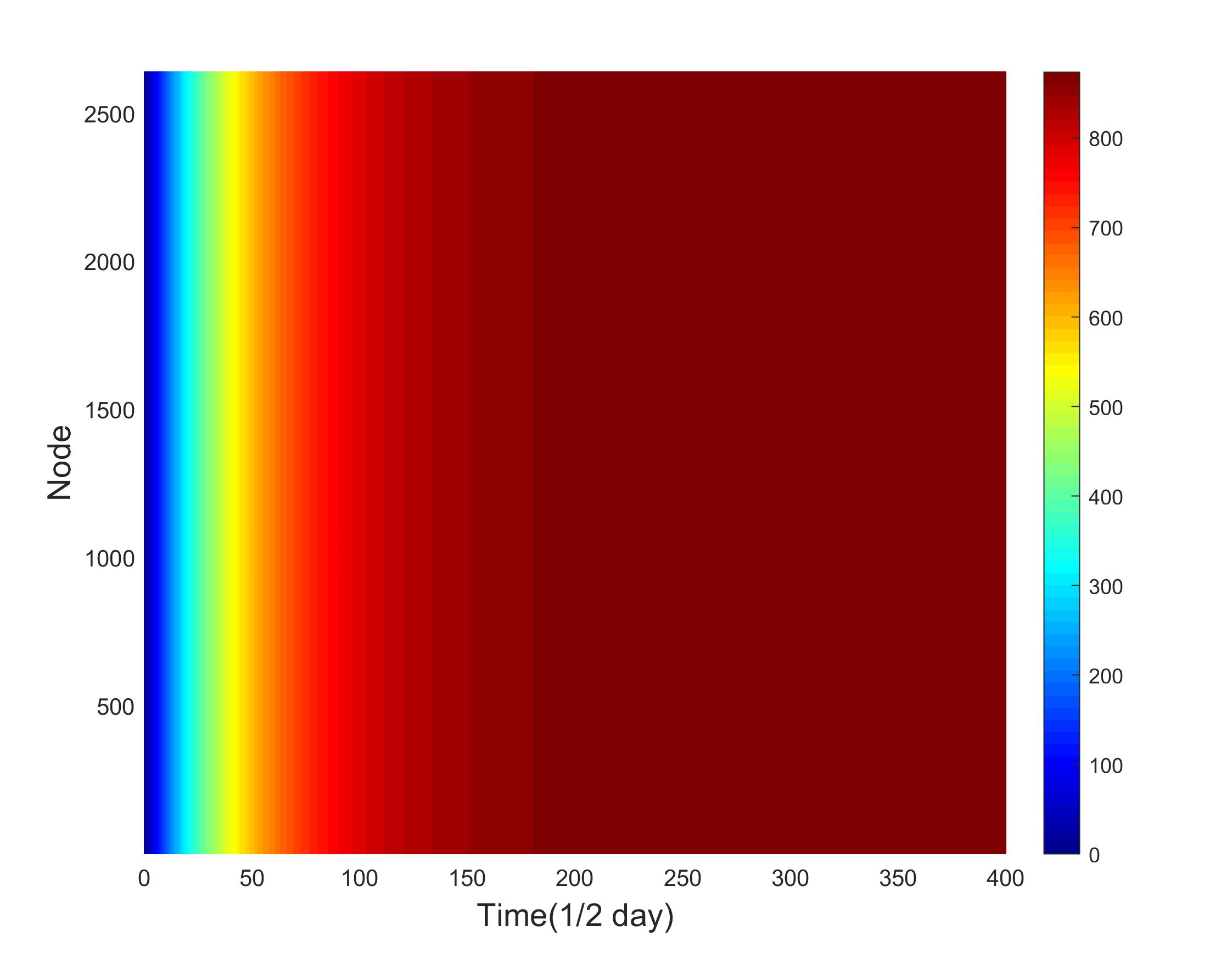}} \\
      \small Solution of $S(x, t)$  &  Solution of $V(x, t)$  \\
      {\includegraphics[width=80mm]{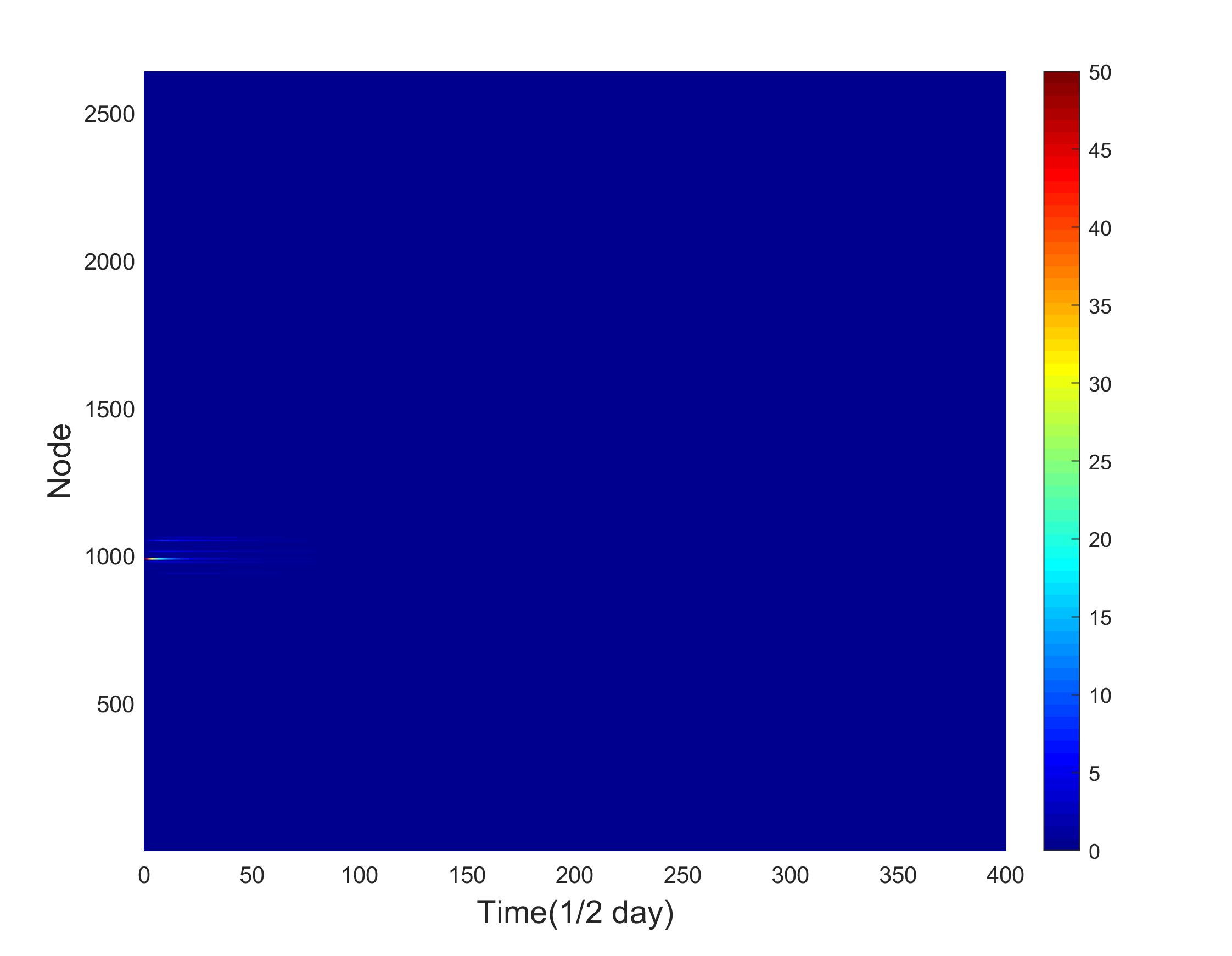}} &
      {\includegraphics[width=80mm]{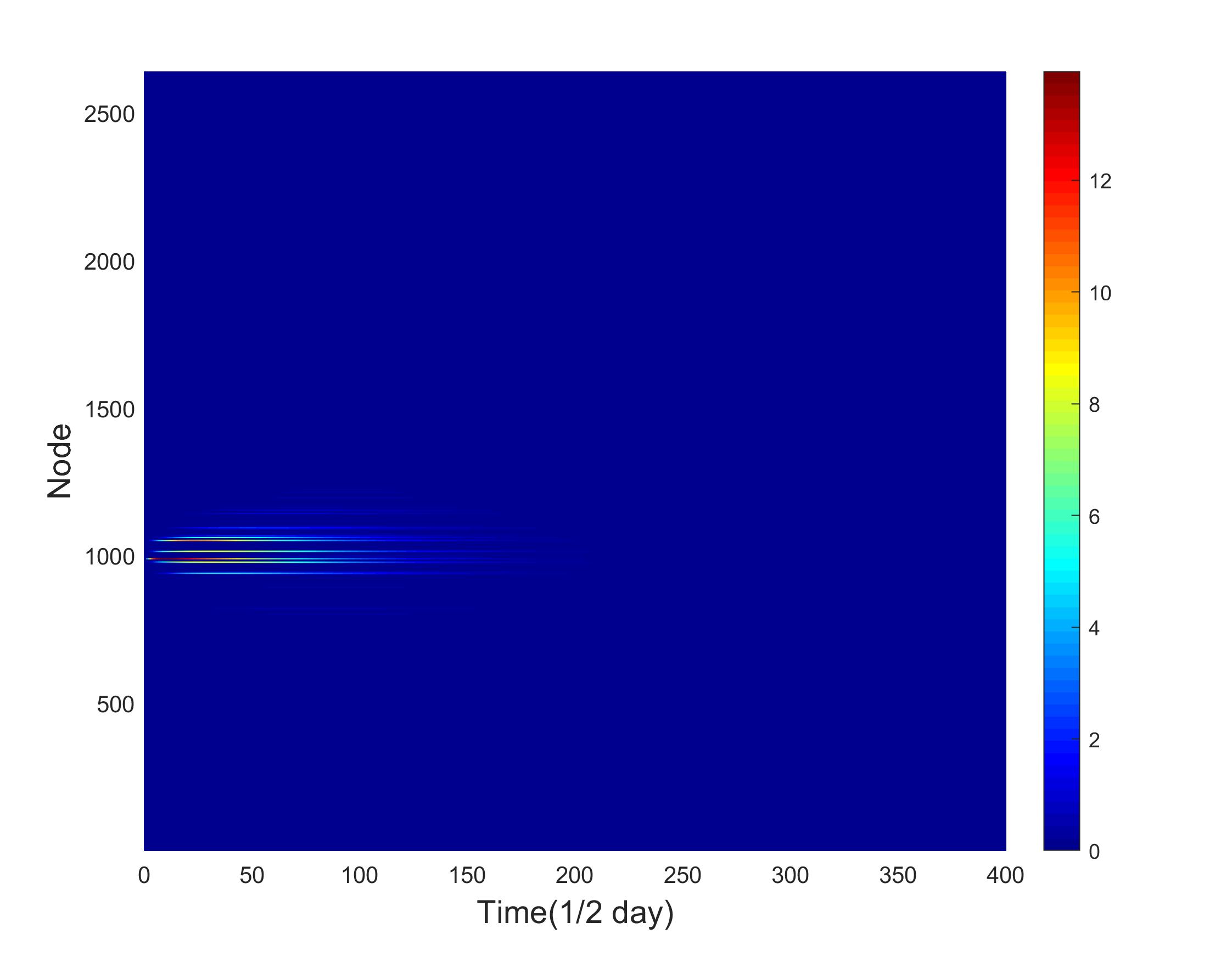}} \\
      \small Solution of $I(x, t)$ & Solution of $R(x, t)$     
\end{tabular}
    \caption{The figure portraits that the solutions of $S(x, t), V(x, t), I(x, t)$ and $R(x, t)$ with $\mathcal{R}_0  = 0.78 < 1$, when $r = 0.01,\, p = 0.002,\, \beta = 0.095/N,\, \xi = 0.02,\, \eta = 0.03,\, \sigma = 0.01,\, \epsilon = 0.05$ and $\gamma = 0.11$.}
    \label{fig:disease_free}
\end{figure}
\begin{figure}[h!]
    \centering
\hspace*{-10pt}\begin{tabular}{c c}
      {\includegraphics[width=80mm]{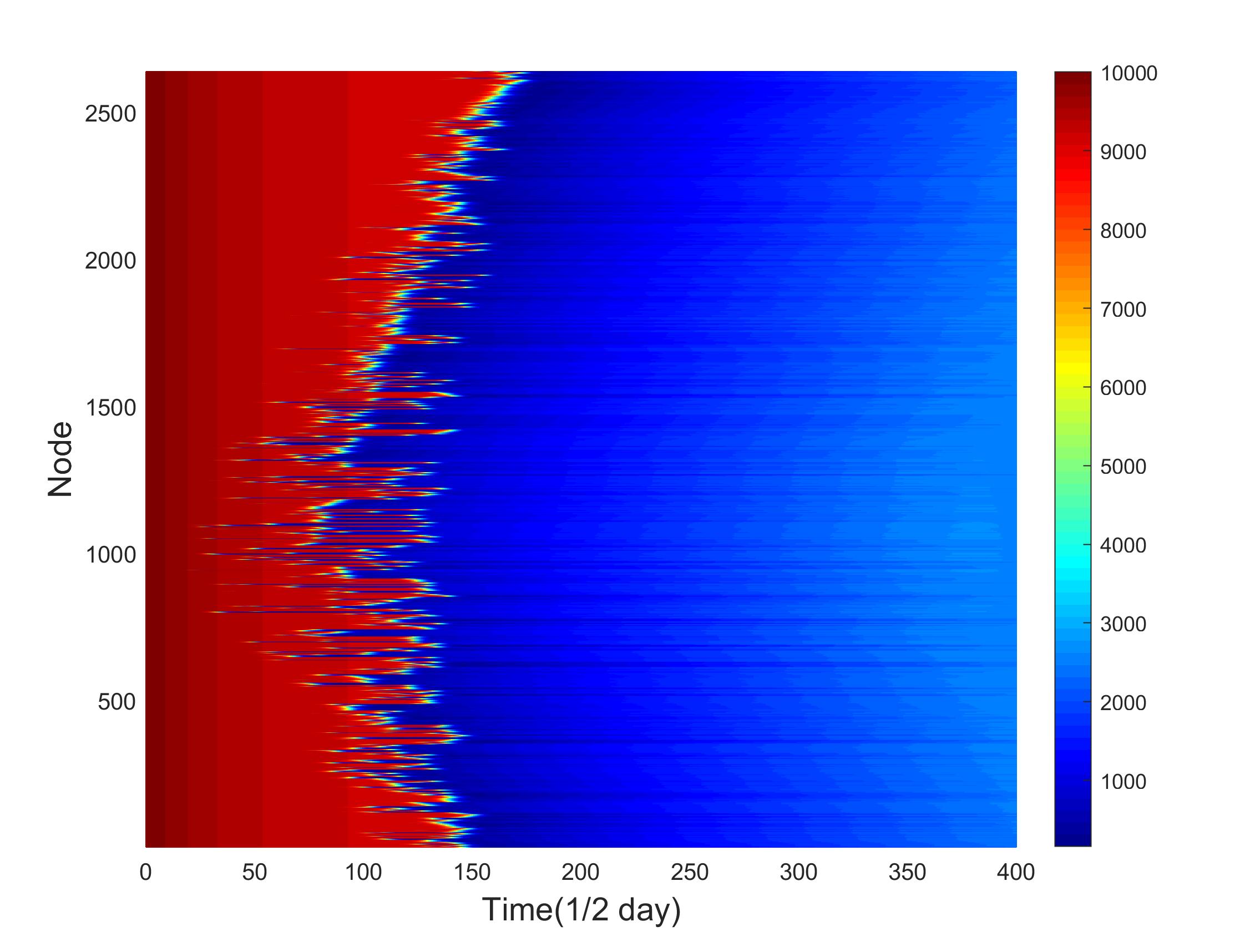}} &
      {\includegraphics[width=80mm]{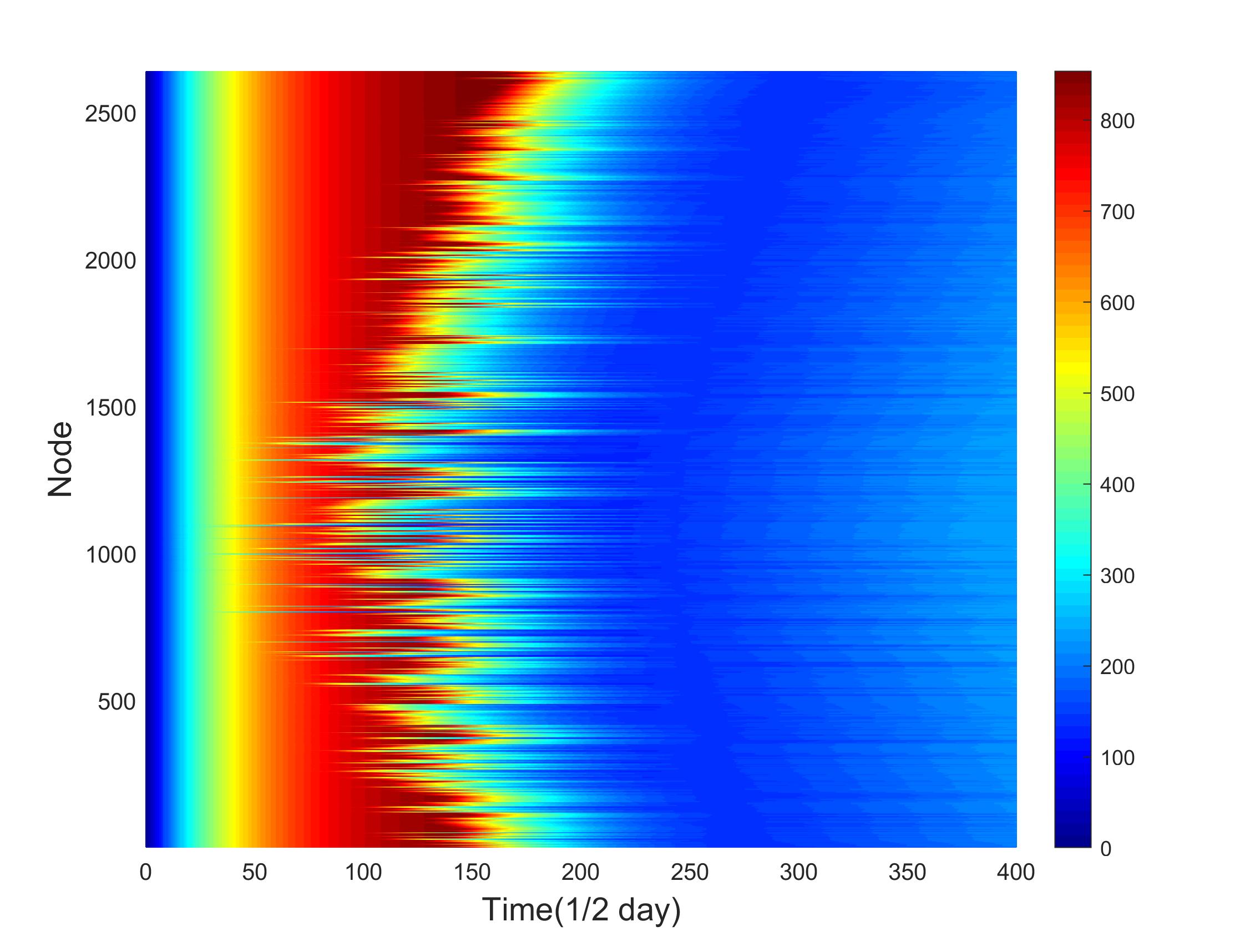}} \\
      \small Solution of $S(x, t)$  &  Solution of $V(x, t)$ \\
      {\includegraphics[width=80mm]{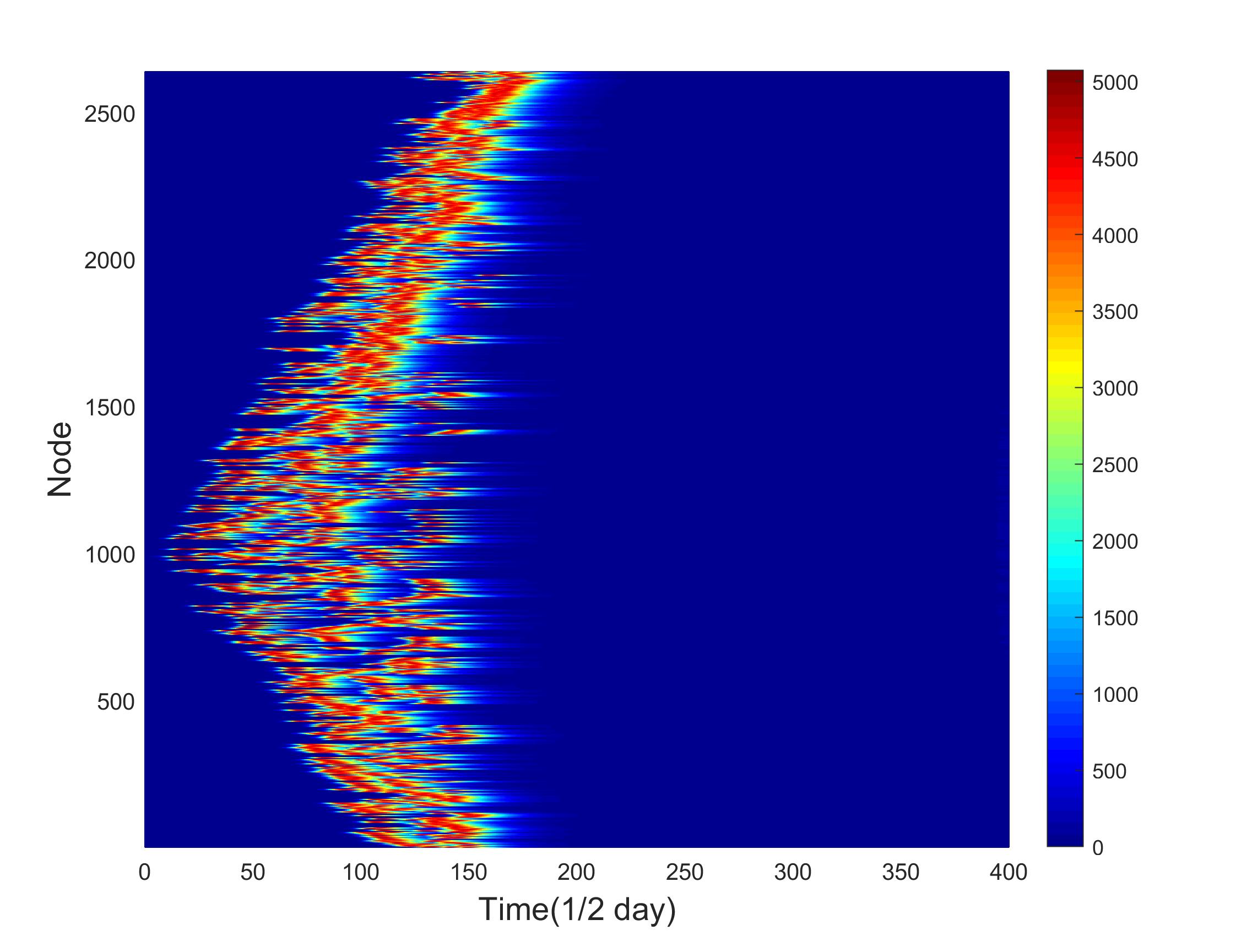}} &
      {\includegraphics[width=80mm]{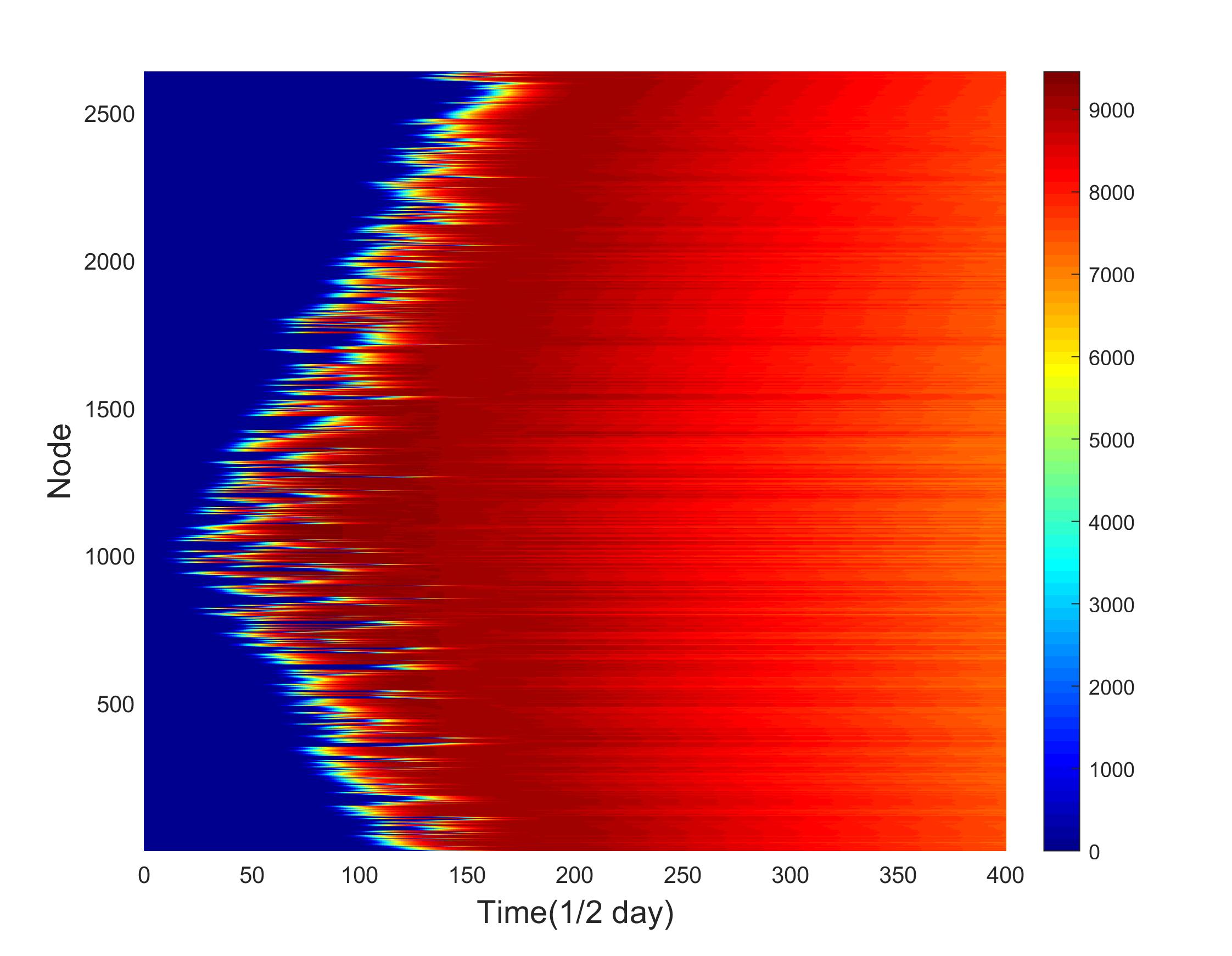}} \\
      \small Solution of $I(x, t)$ & Solution of $R(x, t)$      
\end{tabular}
    \caption{The figure portraits that the solutions of $S(x, t), V(x, t), I(x, t)$ and $R(x, t)$ with $\mathcal{R}_0  = 5.35 > 1$, when $r = 0.01,\, p = 0.002,\, \beta = 0.65/N,\, \xi = 0.02,\, \eta = 0.0,\, \sigma = 0.01,\, \epsilon = 0.05$ and $\gamma = 0.11$.}
    \label{fig:endemic}
\end{figure}
\noindent In Figure \ref{fig:var}, the population mobility rate is fixed $\epsilon = 0.05$. The effect of the parameters $\xi, \sigma, \eta$ are depicted in this figure while the vaccination rate ($p$) is fixed at 0.02 and a fraction of newborns vaccinated is fixed at $r=0.01$. The sub-figure(upper-left) presents the solution of $I(x,t)$ when $\xi, \sigma, \eta$ have no influence in the network. In this case, the disease spreads only a few nodes that are neighbors to the source node, and the maximum peak of infected is around 3,398.  Suppose we only consider that the vaccination is not completely effective and individuals are becoming again infected. Let $\sigma = 0.7$ but $\xi = \eta = 0$, then the disease starts spreading to the entire network with the highest peak of 4,955 within a short period, see sub-figure(upper-right). Next, we only allow the vaccinated individuals are failing to gain immunity against the disease, and they are becoming susceptible individuals again.  Let $\xi = 0.7$ but $\sigma = \eta = 0$, the disease lay out the whole network within a very short time with the highest peak of 5,124, see sub-figure(lower-left). Interestingly, we observe some small waves occur at some periods, but eventually, waves will die. Finally, we assume only that some individuals from the removal class are losing their immunity and becoming again susceptible. Let $\eta = 0.7$ but $\sigma = \xi = 0$, the disease does not spread like the other two parameters $\sigma$ and $\xi$, see sub-figure(lower-right).\\[2ex]
\noindent Fig. \ref{fig:disease_free} demonstrates that all the solutions of the network converge to the disease-free equilibrium point $\textbf{D} = (9126.1, 873.9, 0, 0)$ uniformly as time $t$ goes to infinity when $\mathcal{R}_0  = 0.78< 1$. Thus, the disease-free equilibrium point is globally asymptotically stable, when $\mathcal{R}_0 < 1$. Fig. \ref{fig:endemic} illustrates that all the solutions of the network converge to the endemic equilibrium point $\textbf{E} = (2271.8, 192.0, 4.0, 7532.2)$ as time $t$  uniformly goes to infinity when $\mathcal{R}_0  = 5.35 > 1$. Thus, endemic equilibrium point is globally asymptotically stable, when $\mathcal{R}_0 >1$. 
\section{Conclusions}\label{conclusion}
\noindent A weighted-undirected networked SVIRS model has been proposed with discrete graph Laplacian diffusion. Laplacian diffusion is added for understanding the spread of an infectious disease on the mobility of the population throughout the network. In this network, the probability of an individual travelling from one location to another in any direction is uniform. The proposed model has two equilibria, they are disease-free and endemic. Disease-free equilibrium always occurs while the existence and uniqueness of endemic equilibrium has been shown. Then we presented the stability theory of the model in terms of the basic reproduction number $\mathcal{R}_0$ with the help of Lyapunov function. When the basic reproduction number is below unity, the disease-free equilibrium point is asymptotically globally stable in the network. The endemic equilibrium point is asymptotically globally stable in the network while the basic reproduction number is above unity and $\eta = 0$. Some important numerical results are demonstrated by taking the network of the state of Minnesota.
\section*{Declarations}
\textbf{Conflict of interest} The authors declare that they have no conflict of interest.
\section*{Acknowledgments}
The authors would like to express our gratitude and thank the anonymous referee for their invaluable suggestions and comments, which significantly contributed to the improvement of our original manuscript. MB acknowledges IIITDM Kancheepuram for the research fellowship and facilities.
\vspace*{3pt}              


\end{document}